\documentclass[11pt,reqno]{article}
\usepackage[T1]{fontenc}
\usepackage[utf8]{inputenc}
\input{epsf.sty}
\usepackage{hyperref}
\usepackage[final]{epsfig}
\usepackage{color}
\usepackage{comment}
\usepackage{amsmath, amsthm, amssymb, dsfont, mathrsfs}
\usepackage{enumerate}
\usepackage{tikz}
\newtheorem {thm}{Theorem}[section]
\newtheorem {prop}[thm]{Proposition}
\newtheorem {lem}[thm]{Lemma}
\newtheorem {cor}[thm]{Corollary}
\theoremstyle{definition}

\def\N{{\Bbb N}}

\def\Z{{\Bbb Z}}
\def\R{{\Bbb R}}
\def\P{{\Bbb P}}

\def\E{{\Bbb E}}
\def\one{{\mathds 1}}

\newcommand{\e}{\mathrm e}
\renewcommand{\d}{\mathrm d}

\newcommand\numberthis{\addtocounter{equation}{1}\tag{\theequation}}
\newcommand{\ProofEnde}{\hfill {$\square$}}

\def\0{{\bf 0}}

\def\dist{{\rm dist}}

\def\a{\alpha}

\def\b{\beta}

\def\phi{\varphi}

\def\g{\gamma}

\def\l{\lambda}

\def\k{\kappa}

\def\L{\Lambda}

\def\T{\T}

\def\PP{{\cal P}}

\def\EE{{\cal E}}

\def\Ecal{\mathcal E}

\begin{document}

\title{Exponential moments for planar tessellations}

\author{
Benedikt Jahnel\footnote{
Weierstrass Institute for Applied Analysis and Stochastics, Mohrenstraße 39
10117 Berlin, Germany
\texttt{Benedikt.Jahnel@wias-berlin.de}}
\, and 
Andr\' as T\' obi\' as\footnote{
Berlin Mathematical School, TU Berlin,
Straße des 17.~Juni 136, 10623 Berlin, Germany
\texttt{Tobias@math.tu-berlin.de}}
}

\newcommand{\CC}[1]{{\color{blue} #1}}
\newcommand{\DE}[1]{{\color{red} #1}}
\newcommand{\CK}[1]{{\color{green} #1}}

\maketitle

\begin{abstract}
In this paper we show existence of all exponential moments for the total edge length in a unit disk for a family of planar tessellations based on stationary point processes. Apart from classical tessellations such as the Poisson--Voronoi, Poisson--Delaunay and Poisson line tessellation, we also treat the Johnson--Mehl tessellation, Manhattan grids, nested versions and Palm versions. As part of our proofs, for some planar tessellations, we also derive existence of exponential moments for the number of cells and the number of edges intersecting the unit disk.
\end{abstract}

\bigskip 


\section{Setting and main results}\label{Exp_Setting}
%%%%%%%%%%%%%%%%%%%%%%%%
%Setup Poisson tessellations
%%%%%%%%%%%%%%%%%%%%%%%%
Random tessellations are a classical subject of stochastic geometry with a very wide range of applications for example in the modeling of telecommunication systems, topological optimization of materials and numerical solutions to PDEs. In this paper we focus on {\em random planar tessellations} $S\subset\R^2$ which are derived deterministically from a {\em stationary point process} $X=\{X_i\}_{i\in I}$. The most famous example here is the {\em planar Poisson--Voronoi tessellation}. 

\medskip
%%%%%%%%%%%%%%%%%%%%%%%%
%Setup Poisson tessellations quantities
%%%%%%%%%%%%%%%%%%%%%%%%
Since several decades, research has been performed to understand statistical properties of various characteristics of $S$ such as the degree distribution of its nodes, the distribution of the area or the perimeter of its cells, etc. For the classical examples, where the underlying point process is given by a {\em Poisson point process} (PPP), it is usually possible to derive first and second moments for these characteristics as a function of the intensity $\l$, see~\cite[Table 5.1.1]{OBSC09} and for example~\cite{M89,M94,MS07}. However, to derive complete and tractable descriptions of the whole distribution of these characteristics is often difficult. 

%%%%%%%%%%%%%%%%%%%%%%%%
%Our interest
%%%%%%%%%%%%%%%%%%%%%%%%
\medskip
In this paper we contribute to this line of research by proving existence of all exponential moments for the distribution of the total edge length in a unit disk. More precisely, let $B_r\subset\R^2$ denote the closed centered disk with radius $r>0$ and let $|S \cap A|=\nu_1(S \cap A)$ denote the random total edge length of the tessellation $S\subset\R^2$ in the Lebesgue measurable volume $A\subset\R^2$, where $\nu_1$ denotes the one-dimensional Hausdorff measure. We show for a large class of tessellations that for all $\a\in\R$ we have that
\begin{align}\label{expmomentsfinite} 
\E[\exp(\alpha |S \cap B_1|)] < \infty.
\end{align}

%%%%%%%%%%%%%%%%%%%%%%%%
%Motivation: LD, Cox percolation, SINR percolation
%%%%%%%%%%%%%%%%%%%%%%%%
\medskip
As a motivation, let us mention that the information on the tail behavior of the distribution of $|S \cap B_1|$ provided by~\eqref{expmomentsfinite} is an important ingredient for example in the large deviations analysis of random tessellations. 
%For example one could be interested in the large deviations behavior of the empirical edge measure in 
%\begin{align*}
%\l^{-2}\sum_{i,j\in I}\one\{(X_i,X_j)\in E\}\delta_{(X_i,X_j)\cap B_1}
%\end{align*}
%where $E=E(S)$ denotes the set of edges in $S$.
%
If additionally the tessellation has sufficiently strong mixing properties, namely that there exists $b>0$ such that $|S\cap A|$ and $|S\cap B|$ are stochastically independent for measurable sets $A,B\subset\R^2$ with $\dist(A,B)=\inf\{|x-y|\colon x\in A, y\in B\}>b$, then the cumulant-generating function 
\begin{align*}
\lim_{n\uparrow\infty}n^{-2}\log\E[\exp(-|S\cap B_n|)]
\end{align*}
exists, see~\cite[Lemma 6.1]{HJC18}. This can be used for example to establish the limiting behavior of the percolation probability for the Boolean model with large radii based on Cox point processes where the intensity measure is given by $|S\cap \d x|$, see~\cite{HJC18}. 
Moreover, existence of exponential moments plays a role in establishing percolation in an SINR graph based on Cox point processes in the case of an unbounded integrable path-loss function,
%in giving proper bounds on the acceptable decay of the path-loss function in the context of percolation for the SINR graph based on Cox point processes, 
see~\cite{T18} for details. 

%%%%%%%%%%%%%%%%%%%%%%%%
%Definitions of the tessellations we use: PLT, PVT, PDT, weighted, iterated, Manhattan, Subgraphs
%%%%%%%%%%%%%%%%%%%%%%%%
\subsection{Tessellations}
\begin{figure}
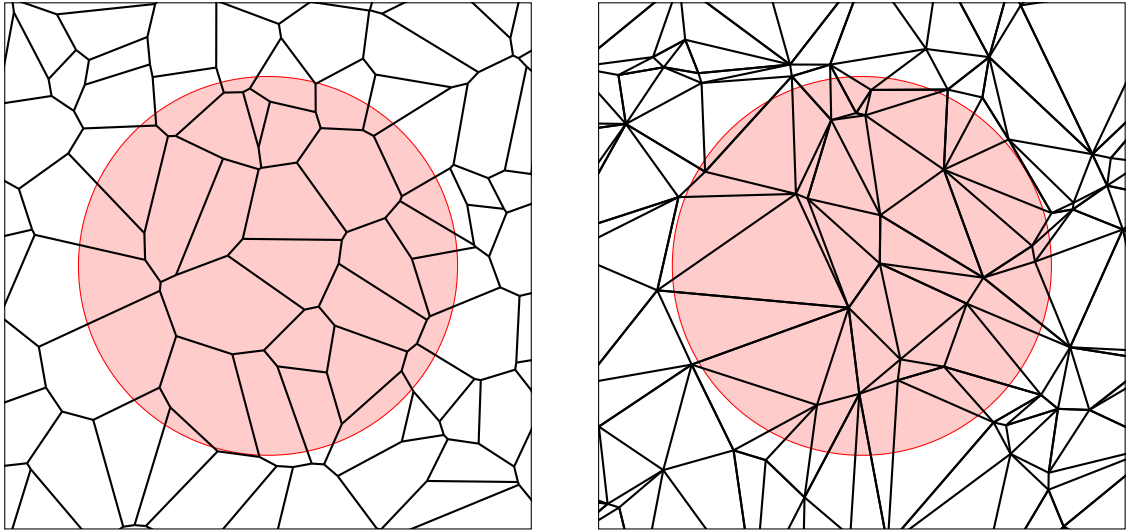

\centering
\input{Fig-PVT.tex}
\hspace{0.5cm}
\input{Fig-PDT.tex}
\caption{A VT (left) and a DT (right) based on a realization of a homogeneous PPP in $\R^2$. The unit disk, in which we estimate the exponential moments of the total edge length, number of edges and number of cells, is shown in red.}
\label{fig-PVTPDT}
\end{figure}

Let $\partial A=\bar A\setminus A^o$ denote the boundary of a set $A\subset \R^2$ and write $x=(x_1,x_2)$ for $x\in\R^2$. Apart from the classical \emph{Voronoi tessellation} (VT), where 
%\begin{align*}
%S_{\rm V}=S_{\rm V}(X)=\{x\in \R^2\colon \text{there exist }X_i\neq X_j\in X\text{ with }|x-X_i|=|x-X_j|\le |x-X_k|\text{ for all }X_k\in X\}
%\end{align*}
\begin{align*}
S_{\rm V}=S_{\rm V}(X)=\bigcup_{i\in I}\partial \{x\in \R^2\colon |x-X_i|=\inf_{j\in I}|x-X_j|\},
\end{align*}
and its dual, the \emph{Delaunay tessellation} (DT), where 
\begin{align*}
S_{\rm D}=S_{\rm D}(X)=\bigcup_{i,j\in I,\,  s\in [0,1]}\{sX_i+(1-s)X_j\colon\exists x\in S_{\rm V}(X)\text{ with }|x-X_i|=|x-X_j|= \inf_{k \in I}|x-X_k|\},
\end{align*}
we also consider the \emph{line tessellation} (LT), where 
\begin{align*}
S_{\rm L}=S_{\rm L}(X)=\bigcup_{i\in I:\, X_i\in \R\times [0,\pi)}\{x\in \R^2\colon x_1\cos X_{i,2}+x_2\sin X_{i,2}=X_{i,1}\}.
\end{align*}
See Figure~\ref{fig-PVTPDT} for realizations of the VT and the DT and their intersections with $B_1$ in case $X$ is a homogeneous PPP.
The extension of the VT known as the {\em Johnson--Mehl tessellation} (JMT) is covered by our results, see for example~\cite{BR08}. For this consider the i.i.d.~marked stationary point process $\widetilde X=\{(X_i,T_i)\}_{i\in I}$ on $\R^2\times [0,\infty)$ with mark measure $\mu(\d t)$. We define the Johnson--Mehl metric by  
\begin{align*}
d_{\rm J}((x,s),(y,t))=|x-y|+|t-s|, \numberthis\label{JMmetric}
\end{align*}
where we use the same notation $|\cdot|$ for the Euclidean norm on $\R^2$ and $[0,\infty)$. 
Then, the JMT is given by 
\begin{align*}
S_{\rm J}&=S_{\rm J}(\widetilde X)=\bigcup_{i\in I}\partial\{x\in\R^2\colon d_{\rm J}((x,0),(X_i,T_i))=\inf_{j \in I} d_{\rm J}((x,0),(X_j,T_j))\}.
%S_{\rm mV}&=S_{\rm mV}(\widetilde X)=\bigcup_{i\in I}\partial\{x\in\R^2\colon d^\phi_{\rm{m}}(x,(X_i,T_i))=\inf_{j \in I} d^\phi_{\rm{m}}(x,(X_j,T_j))\}.
\end{align*}
%For the case $\phi(x)=x$, we for example recover the Johnson--Mehl tessellation, see~\cite{BR08}. For $\phi(x)=x^2$, the associated tessellation is referred to as the Laguerre tessellation, see~\cite{LZ08}.

We also consider the \emph{Manhattan grid} (MG), see for example~\cite{HHJC19}. For this let $Y=(Y_{\rm v},Y_{\rm h})$ be the tuple where $Y_{\rm v}=\{Y_{i,\rm v}\}_{i\in I_{\rm v}}$ and $Y_{\rm h}=\{Y_{i,\rm h}\}_{i\in I_{\rm h}}$ are two independent simple stationary point processes on $\R$. 
Then the MG is defined as
\begin{align*}
S_{\rm M}=S_{\rm M}(Y)=\bigcup_{i\in I_{\rm v},\, j\in I_{\rm h}}(\R\times \{Y_{i,\rm h}\})\cup(\{Y_{j,\rm v}\}\times\R).
\end{align*}
Note that $S_{\rm M}$ is stationary, similarly to all previously defined tessellations, however, unlike them, it is not isotropic. One can make $S_{\rm M}$ isotropic by choosing a uniform random angle in $[0,2\pi)$, independent of $Y$, and rotating $S_{\rm M}$ by this angle. Our results for the MG will be easily seen to hold for both the isotropic and anisotropic version of the MG.

\medskip
Next, let us denote by $(C_i)_{i\in J}$ the collection of cells in the tessellation $S$, where $J=J(S)$. Formally, a cell $C_i$ of $S$ is defined as an open subset of $\R^2$ such that $C_i\cap S=\emptyset$ and $\partial C_i\subset S$.
In view of applications, see for example~\cite{HHJC19, NHGS14}, it is sometimes desirable to consider \emph{nested tessellations} (NT), which we can partially treat with our techniques. For this, let $S_o$ be one of the tessellation processes introduced above, defined via the point process $X^{(o)}$, with cells $(C_i)_{i\in J}$, which now serves as a first-layer process. For every $i\in J$, let $S_i$ be an independent copy of one of the above tessellation processes, maybe of the same type as $S_o$ with potentially different intensity or maybe of a different type, but all $S_i$ should be of the same type and have the same intensity. Let $X^{(i)}$ denote the underlying independent point process of $S_i$. Then the associated NT is defined as
\begin{align*}
S_{\rm N}=S_{\rm N}(X^{(o)},X^{(1)},\dots )=S_o\cup \bigcup_{i\in J}(S_i\cap C_i).
\end{align*}
Here, $\bigcup_{i\in J}(S_i\cap C_i)$ will be called the second-layer tessellation. This definition of a NT originates from \cite[Section 3.4.4]{V09}, where this class of tessellations was defined as a special case of \emph{iterated tessellations}.

\medskip
Note that all kinds of tessellations $S$ defined in this section are stationary, i.e., $S$ equals $S+x$ in distribution for any given $x \in \R^2$. However, for a planar tessellation in order to be stationary, it is not required that it is based on a stationary point process. For example, let $Y$ be a homogeneous PPP on $\R$, and let $X=\{ (X_i,0) \colon X_i \in Y \}$. Then $X$ is not a stationary point process in $\R^2$, however, the associated process $\{ (X_i,t) \colon X_i\in Y, t \in \R \}$ of infinite vertical lines is a stationary tessellation.

\medskip
Finally note that all subgraphs of tessellations having the property~\eqref{expmomentsfinite} inherit this property by monotonicity. In particular, our results cover the cases of the Gabriel graph, the relative neighborhood graph, and the Euclidean minimum spanning tree, since they are subgraphs of the DT, presented in decreasing order with respect to inclusion. 

%%%%%%%%%%%%%%%%%%%%%%%%
%Main theorem (no weighted or interated)
%%%%%%%%%%%%%%%%%%%%%%%%
\subsection{Assumptions}\label{sec-ass}
Unless noted otherwise, throughout the manuscript $X=\{X_i\}_{i\in I}$ denotes a stationary point process on $\R^2$ with intensity $0<\l<\infty$. 
Our results will use the following assumptions on exponential moments for the number of points and void probabilities for the underlying stationary point process. First, for the VT, we assume that
\begin{align}\label{expmomentspointsfinite} 
\limsup_{n\uparrow\infty}|B_{n+4}\setminus B_n|^{-1}\log\E\big[\exp(\beta\#(X\cap B_{n+4}\setminus B_n )\big)\big] < \infty, 
\end{align}
for all $\beta> 0$. Second, we assume that
\begin{align}\label{expmomentsvoidfinite} 
\limsup_{n\uparrow\infty}|B_n|^{-1}\log\P\big(\#(X\cap B_n)=0\big) < 0.
\end{align}
We provide the easy proof that these conditions hold for the homogeneous PPP in Section~\ref{sec-Cox}. They can also be verified for example for some $b$-dependent Cox point processes and some Gibbsian point processes, see also Section~\ref{sec-Cox}.

%\medskip
%For the LT we only have to assume that $\#(X_{1}\cap [0,1])$ has all exponential moments. This is for example the case if $X$ is a PPP. 

\medskip
For the JMT, we generally assume that the mark distribution $\mu(\d t)$ is absolutely continuous with respect to the Lebesgue measure. Further, let $B^{\rm J}_r$ denote the centered ball in the JM metric as defined in~\eqref{JMmetric}. Then, in analogy to the above, we assume that 
\begin{align}\label{expmomentspointsfiniteJM} 
\limsup_{n\uparrow\infty}|B^{\rm J}_{n+4}\setminus B^{\rm J}_n|^{-1}\log\E\big[\exp(\beta\#(\tilde X\cap B^{\rm J}_{n+4}\setminus B^{\rm J}_n )\big)\big] < \infty, 
\end{align}
for all $\beta> 0$. Second, we assume that
\begin{align}\label{expmomentsvoidfiniteJM} 
\limsup_{n\uparrow\infty}|B^{\rm J}_n|^{-1}\log\P\big(\#(\tilde X\cap B^{\rm J}_n)=0\big) < 0.
\end{align}
Again, these conditions hold if $(X_i)_{i \in I}$ is homogeneous PPP and $\mu$ is for example the Lebesgue measure, see Section~\ref{sec-Cox}.

\medskip
For the LT, we will assume that there exists $\beta_\star \leq \infty$ such that the random variable $\# (X \cap ([-1,1] \times [0,2\pi]))$ has exponential moments up to $\beta_\star$, i.e.,
\[ \E\big[\exp\big(\beta \# (X \cap ([-1,1] \times [0,2\pi]))\big)\big] < \infty, \numberthis\label{expmomentsfiniteLT} \]
for all $\beta < \beta_\star$. This condition holds for example for the homogeneous PPP with $\beta_\star=\infty$.

\medskip
For the MG, we assume that there exist $\beta_{\rm v},\beta_{\rm h}\le \infty$ such that the random variables $\#(Y_{\rm v}\cap[0,1])$ and $\#(Y_{\rm h}\cap[0,1])$ have all exponential moments up to $\beta_{\rm v},\beta_{\rm h}$, i.e., 
\begin{align}\label{expmomentsfiniteMG} 
\E[\exp(\beta \#(Y_{\rm v}\cap[0,1]))] < \infty\qquad\text{ and }\qquad\E[\exp(\beta \#(Y_{\rm h}\cap[0,1]))] < \infty,
\end{align}
for all $\beta<\beta_{\rm v}$, respectively $\beta<\beta_{\rm h}$. This condition is satisfied with $\beta_{\rm v}=\infty$ if $Y_{\rm v}$ is a homogeneous Poisson process, and analogously for $\beta_{\rm h}$.

\subsection{Results}
Having defined the types of tessellations we consider, we can now state our main theorem with its proof and all other proofs presented in Section~\ref{Exp_Proofs}.
\begin{thm}\label{thm-expmoments}
We have that~\eqref{expmomentsfinite} holds for all $\a\in\R$ if $S$ is a
\begin{enumerate}
\item Voronoi tessellation, in case~\eqref{expmomentspointsfinite} and \eqref{expmomentsvoidfinite} hold for all $\beta>0$,   \label{thm-PVT}
\item Johnson--Mehl tessellation, in case~\eqref{expmomentspointsfiniteJM} and \eqref{expmomentsvoidfiniteJM} hold for all $\beta>0$,  \label{thm-JMT}
\item Delaunay tessellation, in case the underlying point process is a homogeneous PPP. \label{thm-PDT}
\end{enumerate}
For the line tessellation, in case \eqref{expmomentsfiniteLT} holds for all $\beta<\beta_\star$, then \eqref{expmomentsfinite} holds for all $\alpha<\beta_\star$. \label{thm-PLT}
For the Manhattan grid, in case~\eqref{expmomentsfiniteMG} holds for all $\beta<\beta_{\rm v}$, respectively $\beta<\beta_{\rm h} $, then~\eqref{expmomentsfinite} holds for all $\a<\min\{\beta_{\rm v},\beta_{\rm h}\}$. \label{thm-MG}
\end{thm}

%\begin{thm}\label{thm-expmoments}
%We have that~\eqref{expmomentsfinite} holds for all $\a\in\R$ if $S$ is a
%
%%\vspace{0.1cm}
%\begin{tabular}{ l l }
%  (i) Poisson--Voronoi tessellation,  \label{thm-PVT} & \hspace{1.3cm} (ii) Johnson--Mehl tessellation, \label{thm-JMT} \\
%  (iii) Poisson--Delaunay tessellation, \label{thm-PDT} & \hspace{1.3cm} (iv) Poisson line tessellation, or \label{thm-PLT}  \\
%  (v) Manhattan grid, in case~\eqref{expmomentsfiniteMG} holds for all $\beta>0$. \label{thm-MG} &   \\
%\end{tabular}
%%\begin{enumerate}
%%\item Poisson--Voronoi tessellation, the  \label{thm-PVT}
%%\item Johnson--Mehl tessellation, the \label{thm-JMT}
%%\item Poisson--Delaunay tessellation, the  \label{thm-PDT}
%%\item Poisson line tessellation, or the \label{thm-PLT}
%%\item Manhattan grid. \label{thm-MG}
%%\end{enumerate}
%\end{thm}
%%%%%%%%%%%%%%%%%%%%%%%%
%Cor1 on weighted
%%%%%%%%%%%%%%%%%%%%%%%%
Note that, using H\" older's inequality and stationarity, the statement of Theorem~\ref{thm-expmoments} and all subsequent results remain true if $B_1$ is replaced by any bounded measurable subset of $\R^2$. 

\medskip
Let us briefly comment on the proof of Theorem~\ref{thm-expmoments}. The proof of the parts for the LT and MG is rather  straightforward. As will become clear from the proof, in case of the MG, an application of H\" older's inequality would give the same result without the independence assumption on the point processes $Y_{\mathrm v},Y_{\mathrm h}$, but we lose some of the exponential moments. The cases for the VT, JMT and DT are more involved. However, the statements follow easily if exponential moments for the corresponding number of edges intersecting $B_1$ can be established. More precisely, let $(E_i)_{i\in K}$ denote the collection of edges in the tessellation $S$, where $K=K(S)$, and 
\[ W = \# \{ i \in K \colon E_i\cap B_1 \neq \emptyset \}, \numberthis\label{Wstardef} \]
the number of edges intersecting $B_1$. In the case of tessellations consisting of infinite lines, just as the LT and the MG, one has to be careful with the definition of $W$. Indeed, in this case an edge is to be understood as a maximally linear portion of a cell boundary. Hence, each infinite line of the tessellation contains an infinite number of collinear edges, and $W$ is bounded from below by the number $W_\infty$ of infinite lines having a non-empty intersection with $B_1$. 

Then, for $S$ being a VT or a DT, the edges of $S$ are straight line segments and hence the intersection of each edge with $B_1$ has length at most 2. Similarly, edges of the JMT are either hyperbolic arcs or straight line segments, see \cite[Property AW2, page 126]{OBSC09}.  By convexity, the intersection of any Johnson--Mehl edge with $B_1$ has length at most $|\partial B_1|=2\pi$. Hence, for the VT, DT or JMT, if
\[ \E [\exp(\alpha W)]<\infty \numberthis\label{expmomentsfinitenumberW} \]
holds for all $\alpha>0$, then so does \eqref{expmomentsfinite} for all $\alpha>0$. If \eqref{expmomentsfinitenumberW} holds for some $\alpha>0$, then so does \eqref{expmomentsfinite} for some $\alpha>0$. The following result establishes exponential moments for $W$ and also the simple consequence that 
\begin{align}
\E[\exp(\alpha V)] < \infty,  \label{expmomentsfinitenumberV}
\end{align}
for some $\a>0$, where 
\[ V = \# \{ i \in J \colon C_i\cap B_1 \neq \emptyset \}, \numberthis\label{Vstardef} \]
is the number of cells intersecting $B_1$. 
%
%
%
%For each of these three choices of $S$, it is plausible to start the proof of Theorem~\ref{thm-expmoments} with verifying that exponential moments exist for the number of cells intersecting $B_1$,
%
%
%
%Indeed, 
%We have the following corollary for the weighted tessellations.
%\begin{cor}\label{cor-weighted}
%For $\phi(x)< x^2+c$ strictly increasing on $[0,\infty)$, the aPVT \color{red}and mPVT \color{black} satisfy~\eqref{expmomentsfinite} for all $\a\ge 0$. \color{red} So far this works only for JM. \color{black}
%\end{cor}
\begin{prop}\label{cor-number}
\begin{enumerate}
\item [(i)]
For Voronoi tessellations or Johnson--Mehl tessellations, based on a stationary point process that satisfies~\eqref{expmomentspointsfinite} and~\eqref{expmomentsvoidfinite}, respectively~\eqref{expmomentspointsfiniteJM} and~\eqref{expmomentsvoidfiniteJM}, for all $\beta>0$, \eqref{expmomentsfinitenumberW}
holds for all $\alpha\in\R$. For the Delaunay tessellation based on a homogeneous Poisson point process, \eqref{expmomentsfinitenumberW} holds for some $\a>0$.
\item [(ii)]
For Voronoi tessellations or Johnson--Mehl tessellations, based on a stationary point process that satisfies~\eqref{expmomentspointsfinite} and~\eqref{expmomentsvoidfinite}, respectively~\eqref{expmomentspointsfiniteJM} and~\eqref{expmomentsvoidfiniteJM}, for all $\beta>0$, \eqref{expmomentsfinitenumberV}
holds for all $\alpha\in\R$. For the Delaunay tessellation based on a homogeneous Poisson point process, \eqref{expmomentsfinitenumberV} holds for some $\a>0$.
\end{enumerate}
\end{prop}
As mentioned above, Theorem~\ref{thm-expmoments} parts (i) and (ii) are immediate consequences of Proposition~\ref{cor-number} part (i) for the corresponding tessellations. However, for the case of the DT, as in part (iii) of Theorem~\ref{thm-expmoments}, we cannot use Proposition~\ref{cor-number} since we do not have a statement for all $\alpha>0$. In order to overcome this difficulty, we first estimate small exponential moments of the total number of edges intersecting with $B_a$ for different values of $a>0$ and then use an additional scaling argument to conclude~\eqref{expmomentsfinite} for all $\alpha>0$. Let us also emphasize that for the DT, we establish the above results only in the case in which the underlying point process is a PPP. It is unclear if exponential moments for the number of edges $W$ and number of cells $V$ intersecting with the unit disk exist for the LT and we make no statements about them.
%%%%%%%%%%%%%%%%%%%%%%%%
%Cor2 on iterated
%%%%%%%%%%%%%%%%%%%%%%%%

\medskip
For the NT, existence of exponential moments for $V$ for the first-layer tessellation can be used to verify~\eqref{expmomentsfinite} for $S_{\rm N}$. More precisely, we have the following result. 
%
%
%the number of cells intersecting $B_1$, exhibits all exponential moments. Note that the random number of cells $I_o$ in $S_o$ does not have to coincide with the number of Poisson points used to define $S_o$. Using this we have the following corollary for $S_{\rm N}$. 

\begin{cor}\label{cor-iterated}
Consider the nested tessellation.
\begin{enumerate}[(i)]
\item\label{first-cor-iterated} If for the first-layer tessellation \eqref{expmomentsfinitenumberV} holds for all $\a\in\R$ and for the second-layer tessellation \eqref{expmomentsfinite} holds for all $\alpha\in \R$, then also $S_{\rm N}$ satisfies~\eqref{expmomentsfinite} for all $\a\in\R$.
\item\label{second-cor-iterated} If for the first-layer tessellation \eqref{expmomentsfinitenumberV} holds for some $\a>0$ and for the second-layer tessellation \eqref{expmomentsfinite} holds for some $\alpha>0$, then also $S_{\rm N}$ satisfies~\eqref{expmomentsfinite} for some $\a>0$.
\end{enumerate}
%
%for 
%\begin{enumerate}[(i)]
%\item\label{first-iterated} The nested tessellation satisfies \eqref{expmomentsfinite} for all $\alpha>0$ if the second layer tessellation satisfies~\eqref{expmomentsfinite} for all $\a> 0$ and $S_o$ is such that
%\[ N^* = \# \{ i \in I_o \colon C_i\cap B_1 \neq \emptyset \} \numberthis\label{Nstardef} \]
%has all exponential moments.
%\item\label{second-iterated} \eqref{Nstardef} holds for the PVT, \color{red} the PDT, the PLT \color{black}, and the MG.
%\end{enumerate}
\end{cor}
%In order to use Corollary~\ref{cor-iterated}, we present the following result on the existence of exponential moments for $N_o$. 
%\begin{thm}
%The random variable~\eqref{Nstardef} exhibits all exponential moments for the
%\begin{enumerate}
%\item Poisson--Voronoi tessellation, the  \label{thm2-PVT}\color{red}
%\item Poisson--Delaunay tessellation, the  \label{thm2-PDT}
%\item Poisson line tessellation, or the \label{thm2-PLT}\color{black}
%\item Manhattan grid. \label{thm2-MG}
%\end{enumerate}
%\end{thm}
As we will explain in Section~\ref{sec-nestMG}, the statement of Proposition~\ref{cor-number} is false for the MG based on independent homogeneous Poisson processes on the axes, despite the fact that \eqref{expmomentsfinite} holds in this case according to Theorem~\ref{thm-MG}. However, in the special case where the NT is composed of MGs in both layers and the second-layer MG is based on independent homogeneous Poisson processes, for this $S_{\rm N}$, we still obtain~\eqref{expmomentsfinite} for all $\a\in\R$. This is the content of the following result. 
\begin{prop}\label{prop-iterated-MG}
Consider the nested tessellation and assume that the second-layer tessellation is given by Manhattan grids based on two independent homogeneous Poisson processes and the first-layer tessellation is also a Manhattan grid satisfying~\eqref{expmomentsfinite} for all $\a\in\R$. Then, \eqref{expmomentsfinite} holds for the nested Manhattan grid also for all $\a\in\R$.
\end{prop}

%%%%%%%%%%%%%%%%%%%%%%%%
%Cor3 on Palm
%%%%%%%%%%%%%%%%%%%%%%%%
Let us mention that for the tessellations studied in Theorem~\ref{thm-expmoments}, considering Palm versions of the underlying point process, at least in the case where it is a homogeneous PPP, does not change existence of all exponential moments. We want to be precise here since there are multiple different possibilities to define Palm measures in this context. For the Poisson--VT, Poisson--JMT and Poisson--DT, we denote by $X^*$ the Palm version of the underlying unmarked PPP and denote by $S^*=S(X^*)$ its associated tessellations. For the Poisson--LT we denote by $X^*$ the Palm version of the underlying PPP only with respect to the first coordinate, i.e., $X^* = X \cup \{ (0,\Phi) \}$, where $\Phi$ is a uniform random angle in $[0,\pi)$ that is independent of $X$. Roughly speaking, this corresponds to $S^*_{\rm L}=S_{\rm L}(X^*)$ being distributed as $S_{\rm L}$ when conditioned to have a line crossing the origin $o$ of $\R^2$ with no fixed angle. 
%For the MG we consider only the case where $Y_{\rm v}$ and $Y_{\rm h}$ are homogeneous PPPs on $\R$ with possibly different intensities $\l_{\rm v},\l_{\rm h}>0$. 
The Palm version of the MG is given by 
\begin{align*}
S_{\rm M}^* = (Y_{\mathrm v} \times \R, Y_{\mathrm h}^* \times \R) \mathds 1 \Big\{ U \leq \frac{\lambda_{\mathrm h}}{\lambda_{\mathrm h}+\lambda_{\mathrm v}} \Big\}  +  (Y_{\mathrm v}^* \times \R, Y_{\mathrm h} \times \R)  \mathds 1 \Big\{ U > \frac{\lambda_{\mathrm h}}{\lambda_{\mathrm h}+\lambda_{\mathrm v}} \Big\},
\numberthis\label{PalmMG}
\end{align*}
where $U$ is an independent uniformly distributed random variable on $[0,1]$ and $Y_\mathrm v^*$ and $Y_\mathrm h^*$ denote the Palm versions of $Y_\mathrm v$ and $Y_\mathrm h$, see~\cite[Section III.B]{HHJC19}. We will recall the notion of the Palm version of a general stationary point process in Section~\ref{sec-Palm}.
Palm distributions of NTs can be defined correspondingly, see for example~\cite{HHJC19,V09}. 
\begin{cor}\label{cor-Palm}
Consider all the tessellations $S$ appearing in Theorem~\ref{thm-expmoments}. If the underlying point processes are homogeneous  Poisson point processes, we also have for all $\a\in\R$ that
\begin{align*}
\E[\exp(\alpha |S^* \cap B_1|)] < \infty. \numberthis\label{Palmball}
\end{align*}
\end{cor}
%%%%%%%%%%%%%%%%%%%%%%%%
%Prior work
%%%%%%%%%%%%%%%%%%%%%%%%
To end this section with a short discussion, let us mention that it is a simple consequence of the works~\cite{Z92,C03,H04} that for all $\a\in\R$
\begin{align*}
\E[\exp(\alpha N^*)] < \infty, \numberthis\label{Palmexpmoments}
\end{align*}
where $N^*$ denotes the number of Poisson--Delaunay edges originating from the origin under the Palm distribution for the underlying PPP. The assertion \eqref{Palmexpmoments} seems similar to the one \eqref{Palmball} for the Poisson--DT, however, $S^* \cap B_1$ can contain segments from many edges that are not adjacent to the origin, in particular also from edges both endpoints of which are situated outside $B_1$. It is an interesting open question whether it is possible to provide a simpler proof of the assertion \eqref{Palmball} for all $\alpha\in\R$ or the assertion \eqref{expmomentsfinite} for all $\alpha\in\R$ for the Poisson--DT based on the fact that \eqref{Palmexpmoments} holds for all $\alpha\in\R$.

\medskip
For Corollary~\ref{cor-Palm}, we provide a case-by-case proof. Let us mention that, for the reverse implication with $S^*=S(X^*)$ for a homogeneous PPP $X$ with intensity $\lambda$, using the inversion formula of Palm calculus \cite[Section 9.4]{LP17} and Hölder's inequality, we can derive the following criterion, 
\begin{align*}
\E\Big[\e^{\alpha |S \cap B_1|}\Big]
%= \int \E \Big[ \one\{x\in \mathcal C_o\} \e^{\alpha |S^*\cap B_1(x)|}\Big]\d x  \\
\leq \lambda\int \P(x\in \mathcal C_o)^{1/2}\E \Big[ \e^{2\alpha |S^*\cap B_1(x)|}\Big]^{1/2}\d x= \lambda\int \e^{-\pi |x|^2\lambda/2}\E \Big[ \e^{2\alpha |S^*\cap B_1(x)|}\Big]^{1/2}\d x,  
\end{align*}
where $\mathcal C_o$ is the Voronoi cell of the origin in the VT $S_{\rm V}(X^*)$. 
%
%\color{red} For Corollary~\ref{cor-Palm}, we provide a case-by-case proof. However,  the reverse implication holds more generally. Namely, if $X$ is a PPP with intensity $\lambda$, then for any tessellation $S(X)$, given that the Palm version $S^*(X)=S(X^*)$ has all exponential moments, the same holds for the original tessellation $S(X)$. Indeed, in this case, \eqref{Palmball} clearly also holds with $B_1$ replaced by any compact set, in particular by $B_2$. Let us denote the Voronoi cell of the origin $o$ in the VT $S_{\rm V}(X^*)$ by $\mathcal C_o$. Using the inversion formula of Palm calculus \cite[Section 9.4]{LP17} and Hölder's inequality,
%for $\alpha>0$ we obtain
%\[ \frac{1}{\lambda} \E[\exp(\alpha |S \cap B_1|)] = \E \Big[ \int_{\mathcal C_o} \e^{\alpha |(S^*+x)\cap B_1(x)|} \d x \Big]  \leq \E \big[ |\mathcal C_o| \e^{\alpha |S^* \cap B_2|} \big] \leq \E \big[ |\mathcal C_o|^2]^{1/2} \E \big[ \e^{2\alpha |S^* \cap B_2|} \big]^{1/2}, \]
%Bene: I think the previous line is not correct. 
%\begin{align*}
%\frac{1}{\lambda} \E[\exp(\alpha |S \cap B_1|)] &= \int \E \Big[ \one\{x\in \mathcal C_o\} \e^{\alpha |S^*\cap B_1(x)|}\Big]\d x  \\
%&\leq \int \P(x\in \mathcal C_o)^{1/2}\E \Big[ \e^{2\alpha |S^*\cap B_1(x)|}\Big]^{1/2}\d x \\
%&= \int \e^{-\pi |x|^2\lambda/2}\E \Big[ \e^{2\alpha |S^*\cap B_1(x)|}\Big]^{1/2}\d x
%\end{align*}
%where the right-hand side is finite because $|\mathcal C_o|$ even has exponential moments (cf.~\cite{Z92}).\color{black}

\medskip
Finally, let us comment on possible generalizations of our results to higher dimensions. In at least three dimensions, it is still true that VTs, DTs and JMTs are exponentially stabilizing, i.e., the probability that a point of the underlying point process outside the ball $B_k$ influences the realization of the tessellation inside $B_1$ decays exponentially fast, which is an important argument in our proofs. However, in the planar case, given that the points inside $B_k$ determine the tessellation inside $B_1$, the total number of edges intersecting with $B_1$ can be bounded by constant times the number of points in the region $B_k$ (see e.g.~Section~\ref{sec-PVTproof} for details). This is in general not true in higher dimensions, which yields the main obstacle for generalizing our results. On the other hand, some of our results extend easily to higher dimensions. For example, defining a higher-dimensional analogue of a MG using independent stationary point processes on all coordinate axes and connecting all these points by edges, an analogue of Theorem~\ref{thm-MG} can easily be derived using arguments similar to the ones of Section~\ref{sec-MGproof}. 

\subsection{Examples: Poisson--, Cox-- and Gibbs--Voronoi tessellations}\label{sec-Cox}
It is easy to check that the assumptions listed in Section~\ref{sec-ass}, are satisfied if the underlying point process is a stationary PPP. Indeed, for~\eqref{expmomentspointsfinite} note that by the Laplace transform, for any measurable $B\subset\R^2$
\begin{align*}
%\limsup_{k\uparrow\infty}|B_{k+4}\setminus B_k|^{-1}\log
\E\big[\exp(\beta\#(X\cap B )\big)\big]=\exp\big((\e^{\beta}-1)\l|B|\big).
\end{align*}
Further, for \eqref{expmomentsvoidfinite} note that the void probability for the PPP is given by 
\begin{align*}
\P\big(\#(X\cap B)=0\big)=\exp\big(-\l|B|\big).
\end{align*}
As for the assumptions~\eqref{expmomentspointsfiniteJM} and~\eqref{expmomentsvoidfiniteJM}, the same arguments can be applied. 

\medskip
It is natural to ask under what conditions existence of exponential moments for the total edge length in the unit disk can be guaranteed for tessellations $S(X)$ where $X$ is not a PPP but some different stationary planar point process. As a starting point for future studies, in this section we present examples for the VT based on a stationary Cox point process (CPP) and a stationary Gibbsian point process (GPP) $X$ where our results guarantee the existence of exponential moments.

\subsubsection{Cox--Voronoi tessellations}
A {\em Cox point process} is a PPP with random intensity measure $\L(\d x)$, see for example~\cite{DVJ08} for details. 
We have the following proposition.
\begin{prop}\label{thm-expmoments_Cox}
Consider $S_{\rm V}(X)$ where $X$ is a stationary Cox point process with intensity measure $\L$ satisfying
%\begin{align}\label{Cox_Exp}
%\limsup_{|B| \uparrow \infty} \frac{1}{|B|} \log \E\big[\exp\big(\a \Lambda(B)\big)\big]\le f(\a)
%\end{align}
\begin{align}
\limsup_{n \uparrow \infty}|B_n|^{-1}\log&\E\big[\exp\big(-\Lambda(B_n)\big)\big]<0\qquad\text{ and }\label{Cox_Exp1}\\
\limsup_{n \uparrow \infty}|B_{n+4}\setminus B_n|^{-1}\log&\E\big[\exp\big(\beta (\Lambda(B_{n+4})-\Lambda(B_n))\big)\big]<\infty\label{Cox_Exp2}
\end{align}
for all $\beta>0$. Then, for $S=S_{\rm V}(X)$, \eqref{expmomentsfinite} holds for all $\a\in\R$.
\end{prop}
\begin{proof}[Proof of Proposition~\ref{thm-expmoments_Cox}]
It suffices to verify the assumptions~\eqref{expmomentspointsfinite} for all $\beta>0$ and~\eqref{expmomentsvoidfinite}. For assumption~\eqref{expmomentsvoidfinite}, note that for any measurable $B\subset\R^2$
\begin{align*}
\P\big(\#(X\cap B)=0\big)=\E[\exp(-\L(B))]
\end{align*}
and thus~\eqref{Cox_Exp1} is precisely what we need. The same argument can be applied for assumption~\eqref{expmomentspointsfinite}. 
\end{proof}
The conditions~\eqref{Cox_Exp1} and~\eqref{Cox_Exp2} hold if $\Lambda(Q_1)$ has all exponential moments and $\Lambda$ is $b$-dependent, where we call $\Lambda$ $b$-dependent if for any two measurable sets $A,B\subset \R^2$ such that $\dist(A,B)=\inf_{x \in A,y\in B} |x-y|>b$, the restrictions $\Lambda|_A$ and $\Lambda|_B$ of $\Lambda$ to $A$ respectively $B$ are independent. Indeed, by stationarity of $\Lambda$, it suffices to verify~\eqref{Cox_Exp1} and~\eqref{Cox_Exp2} with $B_k$ replaced by $Q_k=[-k/2,k/2]^2$ (both for $k=n$ and $k=n+4$) in the limit $\N \ni k \to \infty$. Let us assume that $\Lambda$ is $b$-dependent. %Then, for fixed $k$, we can partition $Q_k$ into a bounded number of disjoint subsets such that each of these subsets consists of  (apart from the boundaries) disjoint copies of $Q_1$ and the restrictions of $\Lambda$ to these copies are pairwise independent. 
Then, there exists $b'=b'(b) \in \N$ such that for any $k \in \N$, $Q_k$ can be partitioned into at most $b'$ disjoint subsets such that each of these subsets consists of (apart from the boundaries) disjoint copies of $Q_1$ such that the restrictions of $\Lambda$ to these copies are mutually independent.
Using this independence and the existence of all exponential moments of $\Lambda(Q_1)$, further applying Hölder's inequality for the collection of partition sets,~\eqref{Cox_Exp1} and~\eqref{Cox_Exp2} follow. 
%Let us provide some examples of $\Lambda$ satisfying these properties.
%
%The conditions~\eqref{Cox_Exp1} and~\eqref{Cox_Exp2} hold for all $\a\in\R$ for example for $b$-dependent almost-surely bounded random intensity measures where for some $c>0$ almost surely $\L(B)\le c|B|$. 
A relevant example for a $b$-dependent and even bounded intensity measure is the {\em modulated PPP} where $\Lambda(\d x) = \d x(\lambda_1 \mathds 1 \{ x \in \Xi \} + \lambda_2 \mathds 1 \{ x \in \Xi^{\mathrm c} \} )$, with $\Xi$ being a Poisson--Boolean-model with bounded radii, and $\lambda_1,\lambda_2 \geq 0$, 
see~\cite[Section 5.2.2]{CSKM13}. 
Another example for which conditions~\eqref{Cox_Exp1} and~\eqref{Cox_Exp2} holds, and which is unbounded, is the shot-noise field, see~\cite[Section 5.6]{CSKM13}, 
where $\L(\d x)=\d x\sum_{i\in I}\k(x-Y_i)$ for some integrable kernel $\k\colon \R^2\to[0,\infty)$ with compact support and $\{Y_i\}_{i\in I}$ a stationary PPP. 
%%%%%%shot-noise field computation of the limit
%Indeed, for the shot-noise field, $\L(B)=\sum_{i\in I}\int\d x \k(x)\one\{x\in B(Y_i)\}\le\#(Y_i\in C) \int\d x\k(x)$ with $C=C'\oplus B$ where $C'$ denotes the support of $\k$ . Now, $Z=\#(Y_i\in C)$ is a Poisson random variable and we denote its parameter $\rho$. Then, 
%\begin{align*}
%\limsup_{|B|\uparrow\infty}|B|^{-1}\log \E[\exp(\a \Lambda(B))]\le\rho(\exp(\a\int\d x\k(x))-1)\limsup_{|B|\uparrow\infty} \frac{|C|}{|B|}=\rho(\exp(\a\int\d x\k(x))-1)
%\end{align*}
%has the desired property. 
\subsubsection{Gibbs--Voronoi tessellations}
A {\em Gibbs point process} on $\R^2$ is defined via its conditional probabilities in bounded measurable volumes $B\subset \R^2$. They take the form of a {\em Boltzmann weight}
\begin{align*}
\PP_B(\d X_B)\frac{\exp\big(-\gamma H(X_BX_{B^{\rm c}})\big)}{\int\PP(\d X'_B)\exp\big(-\gamma H(X'_BX_{B^{\rm c}})\big)},
\end{align*}
where $\PP_B$ is a PPP on $B$ with intensity $\l>0$, $\gamma\in\R$ is a system parameter and $H$ is the {\em Hamiltonian}, which assigns some real-valued energy to the configuration $X_BX_{B^{\rm c}}=X_B\cup X_{B^{\rm c}}$, where $X_{B^{\rm c}}$ is a boundary configuration in $B^{\rm c}=\R^2\setminus B$. For details see for instance~\cite{D19}. As an example, we consider the {\em Widom--Rowlinson model} where $H(X)=|\bigcup_{X_i\in X}B_r(X_i)|$, with $B_r(x)$ the ball of radius $r>0$, centered at $x\in\R^2$. Existence of associated point processes on $\R^2$ that are stationary can be guaranteed, see for example~\cite{CCK95}. We have the following result. 

\begin{prop}\label{thm-expmoments_Gibbs}
Consider $S_{\rm V}(X)$ where $X$ is the Widom--Rowlinson model. Then, for $S=S_{\rm V}(X)$, \eqref{expmomentsfinite} holds for all $\a\in\R$.
\end{prop}
\begin{proof}[Proof of Proposition~\ref{thm-expmoments_Gibbs}]
It suffices to verify the assumptions~\eqref{expmomentspointsfinite} for all $\beta>0$ and~\eqref{expmomentsvoidfinite}. For assumption~\eqref{expmomentsvoidfinite}, note that by consistency for all bounded measurable $B\subset\R^2$, 
\begin{align*}
|B|^{-1}\log\P\big(\#(X\cap B)=0\big)&=|B|^{-1}\log\frac{\exp(-\l|B|)\exp(-\g H(X_{B^{\rm c}}))}{\int\PP(\d X'_B)\exp(-\g H(X'_BX_{B^{\rm c}}))}\\
&\le -|B|^{-1}\log\sum_{n\ge 0}\tfrac{1}{n!}(\l|B|)^n\exp(-n\gamma\pi r^2)=-\l\exp(-\gamma\pi r^2)<0,
\end{align*}
for all $r,\l,\gamma>0$. For assumption~\eqref{expmomentspointsfinite}, note that for all bounded measurable sets $B\subset\R^2$,
\begin{align*}
|B|^{-1}\log\E\big[\exp(\beta\#(X\cap B))\big]\le|B|^{-1}\log\frac{\sum_{n\ge 0}\tfrac{1}{n!}(\l|B|)^n\exp(\beta n)}{\sum_{n\ge 0}\tfrac{1}{n!}(\l|B|)^n\exp(-n\gamma\pi r^2)}=\l(\e^{\beta}-\e^{-\gamma\pi r^2})<\infty,
\end{align*}
for all $r,\l,\gamma,\beta>0$, which proves the desired result. 
\end{proof}

\subsection{Absence of exponential moments for the number of edges and cells}\label{sec-nestMG}
In Proposition~\ref{cor-number}, we provide statements about existence of exponential moments for $V$, the number of cells intersecting $B_1$, and $W$, the number of edges intersecting $B_1$. In this section we want to exhibit one example in our family of tessellations for which exponential moments for $V$ do not exist. Indeed, take the MG where the underlying stationary point processes are independent PPPs $Y_{\rm v}$ and $Y_{\rm h}$ with intensity $\l$. By translation invariance, we can also consider the random variable $V'$, the number of cells intersecting $Q_1$.  In order to simplify the notation, let us write $X_\mathrm v=\#(Y_{\rm v} \cap [-1/2,1/2])$ and $X_\mathrm h=\#(Y_{\rm h} \cap [-1/2,1/2])$. These random variables are independent and Poisson distributed with parameter $\lambda$. Then we have that 
\begin{align*}
\E[\exp(\a V')]&=\e^{\alpha} \E\big[\exp\big(\alpha ( X_\mathrm v +X_\mathrm h + X_{\rm v} X_{\rm h} )\big)\big] =\e^{\alpha} \sum_{k=0}^{\infty}\E \big[ \exp\big(\alpha ( k +(k+1) X_\mathrm h )\big)\big] \P(X_{\rm v}=k)\\
& = \e^{\alpha} \sum_{k=0}^{\infty} \e^{\alpha k} \frac{\lambda^k}{k!} \e^{-\lambda} \exp\big(\lambda (\e^{\alpha(k+1)}-1)\big) = \e^{\alpha-2\lambda} \sum_{k=0}^{\infty} \exp \big(\alpha k+\lambda \e^{\alpha(k+1)} \big) \frac{\lambda^k}{k!}=\infty.
\end{align*}
Since for the MG based on PPPs, $V$ and the number $W_\infty$ of infinite lines intersecting with $Q_1$ are of the same order, further, $W_\infty \leq W$, it follows that $\E[\exp(\a W)]=\infty$.

\section{Proofs}\label{Exp_Proofs}
%We will use the following notations in this section. We let $Q_r$, respectively $B_r$, denote the box of side length $r>0$, respectively the closed ball of radius $r>0$, centered at the origin $o$. For $A \subseteq \R$, we will write $\overline A$ for the closure of $A$.
For our results, it obviously suffices to consider $\a>0$ instead of $\a\in\R$. 

%%%%%%%%%%%%%%%%%%%%%%%%
%PLT, PVT, PDT
%%%%%%%%%%%%%%%%%%%%%%%%

\subsection{Total edge length, number of edges and cells: Proof of Theorem~\ref{thm-expmoments} and Proposition~\ref{cor-number}}
The proof of Theorem~\ref{thm-expmoments} and Proposition~\ref{cor-number} is organized as follows. As already discussed, for the VT and the JMT it suffices to show Proposition~\ref{cor-number} part (i) for all $\alpha>0$ in order to conclude the corresponding part of Theorem~\ref{thm-expmoments}. 
In Section~\ref{sec-PVTproof} we cary out the proof of Proposition~\ref{cor-number} part (i) for all $\alpha>0$ for the VT and in \ref{sec-JMTproof} for the JMT. Section~\ref{sec-PDTproof} is devoted to the case of the Poisson--DT. Here we first verify an extended version of Proposition~\ref{cor-number} part (i) for small $\alpha>0$, and using this we verify \eqref{expmomentsfinite} for all $\alpha>0$. The direct and short proofs of \eqref{expmomentsfinite} for all $\alpha>0$ for the LT and the MG can be found in Sections~\ref{sec-PLTproof} and \ref{sec-MGproof}, respectively. Given these results, we prove Proposition~\ref{cor-number} part (ii) in Section~\ref{sec-cells}.
\subsubsection{Voronoi tessellations: Proof of Proposition~\ref{cor-number} part (i)}\label{sec-PVTproof}
%Let $B_r$ denote the centered ball with radius $r>0$ and $|S^\l\cap A|=\nu_1(S^\l\cap A)$, the random edge length of the Poisson Voronoi tessellation $S$ in the volume $A\subset\R^2$, where $\nu_1$ denotes the one-dimensional Hausdorff measure. 
%\begin{thm}\label{thm-Voronoiexpmoments}
%For all $\a>0$, $\E[\exp(\a|S \cap B_1|)]<\infty$.
%\end{thm}
%\begin{proof}
%We follow the same arguments as in AT18. The main new idea is to estimate numbers of Poisson points in balls extended only by a constant difference of radii. 
It suffices to verify \eqref{expmomentsfinitenumberW} for all $\alpha>0$. For this we extend arguments first presented in \cite[Theorem 2.6]{T18}, where it was shown that $\mathbb E[\exp(\alpha |S_{\rm V} \cap [-1/2,1/2]^2|)]<\infty$ for some $\alpha>0$ in the case where the underlying point process is a PPP. 
Let us extend the notion of $W$ defined in \eqref{Wstardef} to balls of different radii via
\[ W_a = \# \{ i \in K \colon E_i\cap B_a \neq \emptyset \}, \numberthis\label{Wadef} \]
where we recall that $(E_i)_{i \in K}$ is the collection of edges of $S_{\rm V}$. 
The following lemma states that unless we have a void space, numbers of edges can be bounded from above by numbers of points in bounded regions. 
\begin{lem}\label{lem-conditioningonapoint}
Let $b\geq a > 0$. If $X\cap B_b \neq \emptyset$, then we have
\[ W_a \leq 3 \# (X \cap B_{b+3a}). \numberthis\label{perimeterbound} \]
\end{lem}
\begin{proof}
Let us assume existence of $X_i\in X\cap B_b$.  We first claim that for any edge of $S_{\rm V}$ intersecting with $B_a$, the corresponding edge in the dual DT connects two points in $B_{b+3a}$. Indeed, assume otherwise, then there exists $v \in B_a$ and $X_j \in X \cap B_{b+3a}^{\mathrm c}$ such that $|v-X_j|=\min \{ |v-X_l| \colon l \in I \}$ and 
\[ |v-X_j| \geq \dist(\{X_j\},B_{a}) > (b+3a)-a >b+a. \]
On the other hand,
\[ |v-X_i| \leq \max_{y \in B_a,z\in B_b} |y-z| = 2 a+(b-a)=b+a, \] 
which is a contradiction.
Thus, for any Voronoi edge intersecting with $B_a \subseteq B_b$, the corresponding Delaunay edge has both endpoints in $X \cap B_{b+3a}$. But since the subgraph of the Delaunay graph spanned by the vertex set $X \cap B_{b+3a}$ is simple and planar, Euler's formula (see e.g.~\cite[Remark 2.1.4]{M94}) implies that the number of such edges is bounded by 3 times the number of vertices in this subgraph. This implies \eqref{perimeterbound}.
\end{proof}

Note that Lemma~\ref{lem-conditioningonapoint} holds for any point cloud $X$. 
The proof of \eqref{expmomentsfinitenumberW} for the VT now rests on the assumption that it is exponentially unlikely to have large void spaces of order $k^2$ and existence of exponential moments for numbers of points in annuli of order $k$.
Let 
\[ R= \inf \{ r>0 \colon B_r \cap X \neq \emptyset \} \numberthis\label{Rdef} \] denote the distance of the closest point in $X$ to the origin. 
%Then we have for all $r>0$ that 
%\[ \mathbb P(R \geq r) = \exp(-\lambda r^2 \pi). \numberthis\label{exponentialstabilization} \]
%this is a version of the statement that the Poisson--Voronoi tessellation is exponentially stabilizing \cite[Example 3.1]{HJC17}. 
\begin{proof}[Proof of Proposition~\ref{cor-number} part (i)] In the event 
$\{ R \leq 1 \}$ we have that $B_1 \cap X \neq \emptyset$, and therefore by Lemma~\ref{lem-conditioningonapoint} applied for $a=b=1$, we obtain
\[ W \leq 3\# (X \cap B_{4}). \]
%The random variable $\# (X \cap B_{4})$ is Poisson distributed with parameter $16\pi\lambda$ and hence has all exponential moments. 
On the other hand, in the event $\{ R \geq 1 \}$, we can apply Lemma~\ref{lem-conditioningonapoint} with $a=1$ and $b=R$ in order to obtain that, almost surely,
\[ W \leq 3 \# (X \cap B_{R+3}) = 3+ 3\# (X \cap (B_{R+3} \setminus B_R)),  \]
where we also used that by stationarity, on $\partial B_R$ there is precisely one point, almost surely.  
By assumption~\eqref{expmomentsvoidfinite} we have $\E[R]<\infty$ and hence $\mathbb P(R<\infty)=1$. 
%and the distribution of $R$ has a density with respect to the Lebesgue measure given by $2\lambda\pi x\e^{-\lambda\pi x^2}$, 
We can thus estimate for all $\a>0$,
%\begin{equation} \label{Poissonestimate}
%\begin{split}
%\E&[\exp(\a W)]
%%=2\lambda\pi\int_0^\infty x\e^{-\lambda\pi x^2}\E[\exp(\a W)|R=x]\d x\cr
%\le \E[\exp(3\a\# (X \cap B_{4}))]+\e^{3\a} 2\lambda\pi \int_1^\infty x\e^{-\lambda\pi x^2}\E\big[\exp\big(3\a \# (X \cap (B_{x+3} \setminus B_x)\big)\big] \d x\cr
%&= \E[\exp(3\a\# (X \cap B_{4}))]+\e^{3\a} 2\lambda\pi \int_1^\infty x\exp\big((6x+9)\l\pi (\exp(3\alpha)-1)-\lambda\pi x^2\big) \d x,
%\end{split}
%\end{equation}
\begin{equation} \label{Poissonestimate}
\begin{split}
&\E[\exp(\a W)]\le \E\big[\exp(3\a\# (X \cap B_{4}))\big]+\e^{3\a}\sum_{k\ge 2}\E\big[\exp\big(3\a\# (X \cap (B_{R+3} \setminus B_R))\big)\one\{R\in[k-1,k)\}\big]\cr
&\le\E\big[\exp(3\a\# (X \cap B_{4}))\big]+\e^{3\a}\sum_{k\ge 2}\E\big[\exp\big(3\a\# (X \cap (B_{k+3} \setminus B_{k-1}))\big)\one\{\# (X \cap B_{k-1})=0\}\big]\cr
&\le\E\big[\exp(3\a\# (X \cap B_{4}))\big]+\e^{3\a}\sum_{k\ge 1}\E\big[\exp\big(6\a\# (X \cap (B_{k+4} \setminus B_{k}))\big)\big]^{1/2}\P\big(\# (X \cap B_{k})=0\big)^{1/2},
%&= \E\big[\exp(3\a\# (X \cap B_{4}))\big]+\e^{3\a} 2\sum_{k\ge 2}\exp\Big(\tfrac{1}{2}\lambda\pi\big((6k+8)(\exp(6\alpha)-1)-(k-1)^2\big)\Big),
\end{split}
\end{equation}
where we used H\" older's inequality in the last line. Now, by the assumptions~\eqref{expmomentsvoidfinite} and~\eqref{expmomentspointsfinite}, there exist $c_1,c_2>0$ such that for sufficiently large $k$, we have 
\begin{equation*}
\begin{split}
\E\big[\exp\big(6\a\# (X \cap (B_{k+4} \setminus B_{k}))\big)\big]\P\big(\# (X \cap B_{k})=0\big)\le \exp\big(\pi(c_1k-c_2k^2)\big),
\end{split}
\end{equation*}
and hence summability of the right-hand side of~\eqref{Poissonestimate} is guaranteed. 
%Hence it suffices to verify that for all $c>0$, we have that $\E [ \exp ( c R ) ]<\infty$. But this is true since, 
%
%
%Now, conditional on $R$, $\# (X \cap (B_{R+3} \setminus B_R))$ is Poisson distributed with parameter $\lambda ((R+3)^2-R^2)\pi=(6R+9)\lambda \pi$, and hence, using the Laplace transform for Poisson random variables, for $\alpha>0$ we have
%\begin{equation} \label{Poissonestimate}
%\begin{split}
%\E[\exp(3\alpha\#(X \cap( B_{R+3} \setminus B_R))) \mathds 1 \{ R > 1 \}]\le\E[\exp\big((6R+9)\l\pi (\exp(3\alpha)-1)\big)].
%%\E\big[\exp(\alpha |S_{\rm V} \cap B_1|) \mathds 1 \{ R>1\} \big] &\leq \exp(2\pi\alpha) 
%%\E \big[\exp(2\pi\alpha Y((2R+6) \l)) \mathds 1 \{ R>1\} \big]  \\
%%&\leq \exp(2\pi\alpha)\E[\exp\big((2R+6)\l (\exp(2\pi\alpha)-1)\big)].
%\end{split}
%\end{equation}
%Finally, note that  for any non-negative random variable $Z$ we have $\E(Z) \leq \sum_{k=0}^{\infty} \P(Z \geq k)$ and thus, \color{red}for $c=3\l (\exp(3\alpha)-1)$\color{black}, \eqref{exponentialstabilization} implies that 
%\[  \mathbb E [ \exp ( c R ) ] \leq 1 + \sum_{k=1}^{\infty} \mathbb P \big( \exp(c R) \geq k \big) =1 + \sum_{k=1}^{\infty} \mathbb P \Big(R \geq \frac{\log k}{c}\Big) =1+ \sum_{k=1}^{\infty} \exp\Big( -\frac{\l \pi (\log k)^2}{c^2}\Big)<\infty, \]
%which can be seen for example by the Leibniz criterion. 
This concludes the proof of Proposition~\ref{cor-number} for the number of edges and thus of Theorem~\ref{thm-expmoments} part~(i). 
 \end{proof}

\subsubsection{Johnson--Mehl tessellations: Proof of Proposition~\ref{cor-number} part (i)}\label{sec-JMTproof}
As explained in Section~\ref{Exp_Setting}, Theorem~\ref{thm-expmoments} (ii) follows once we verify \eqref{expmomentsfinitenumberW} for the JMT for all $\alpha>0$, which is the first part of Proposition~\ref{cor-number} for the JMT. 
%We fix $\varphi$ satisfying the condition of the corollary (i.e., $\phi$ strictly increasing with $\phi(x)<x^2+c$). Let us define a metric $d$ on $\R^3=\R^2 \times\R$ as follows
%\[ d((x,t),(y,u))=\varphi(|x-y|) + \varphi(|t-u|), \qquad (x,t),(y,u)\in\R^3. \]
%Note that if $t \geq 0$, we have that $d_{\text a}^\varphi(y,(x,t))=d((y,u),(x,t))$. For $(x,t) \in \R^3$ and $r>0$, we write $B_r^{d}(x,t)$ for the open ball of radius $r$ centered at $(x,t)$ in $\R^3$ equipped with the $d$-distance. Further, we write
The arguments are very similar to the ones used in Section~\ref{sec-PVTproof} for the VT. To start with, we have the following lemma, which is an analogue of Lemma~\ref{lem-conditioningonapoint} in the Johnson--Mehl case. Recall that for $(x,s) \in \R^2 \times [0,\infty)$ and $r>0$ we write $B_r^{\rm J}(x,s)$ for the closed ball of radius $r$ around $(x,s)$ in the Johnson--Mehl metric, see~\eqref{JMmetric}.

\begin{lem}\label{lem-JMindependence}
Let $b \geq a>0$. If $ \widetilde X \cap B^{\rm J}_b \neq \emptyset$, then $S_{\rm J} \cap B_a$ is determined by $\widetilde X \cap B_{b+3a}^{\mathrm J}$. That is, for any $x\in S_{\rm J} \cap B_a$, if $j \in I$ is such that $d_{\rm J}((X_j,T_j),(x,0))=\inf_{k\in I} d_{\rm J} ((X_k,T_k),(x,0))$, then $(X_j,T_j) \in B_{b+3a}^{\mathrm J}$. 
\end{lem}
\begin{proof}
Assume that there exists $i \in I$ such that $(X_i,T_i) \in B^{\rm J}_b$ and that $S_{\rm J}$ exhibits an edge having a non-empty intersection with $B_a$, and let $x \in B_a$ be a point of such an edge. Then, using the triangle inequality, since
\[ d_{\rm J}((x,0),(X_i,T_i)) \leq d_{\rm J}((x,0),(o,0))+d_{\rm J}((o,0),(X_i,T_i)) =|x|+d_{\rm J}((o,0),(X_i,T_i))\leq a+b, \]
and for any $j \in I$ with $(X_j,T_j) \notin B^{\rm J}_{b+3a}$, we have
\begin{align*}
 d_{\rm J}((x,0),(X_j,T_j)) & = T_j + |x-X_j| \geq (T_j + |X_j|)-|x| > b+3a-a=b+2a>b+a \\
 & \geq d_{\rm J}((x,0),(X_i,T_i)), 
\end{align*}
and the result follows.
\end{proof}
\begin{proof}[Proof of \eqref{expmomentsfinitenumberW} for the JMT for all $\alpha>0$]

We start with two preliminary observations. First, let $\Ecal$ denote the set of (closed) edges of $S_{\rm J}$. 
%, so that $S_{\mathrm J}=(\Vcal,\Ecal)$. 
By construction of a JMT, almost surely, any $E \in \Ecal$ has the property that there exist precisely two points $(X_i,T_i),(X_j,T_j)$ (depending on $E$) such that for all $z \in E$
\[ d_{\mathrm J}((z,0),(X_i,T_i))=d_{\mathrm J}((z,0),(X_j,T_j))=\inf_{k \in I} d_{\mathrm J}((z,0),(X_k,T_k)). \]
In this case, we will write $E=\big((X_i,T_i);(X_j,T_j)\big)$. We claim that for any finite subset $I_0$ of $I$, 
\[ \# \{ \big((X_i,T_i);(X_j,T_j)\big) \in \Ecal \colon i,j \in I_0 \} \leq 3 \# I_0 \numberthis\label{dualityargument}  \]
holds. Indeed, the set on the left-hand side of \eqref{dualityargument} is in one-to-one correspondency with $\#\mathcal D(I_0) $ where 
 \[ \mathcal D(I_0) =\{ (i,j) \in I^2_0 \colon \big((X_i,T_i);(X_j,T_j)\big) \in \EE \}, \]
 since $(i,j) \in \mathcal D(I_0)$ if and only if $X_i$ and $X_j$ are connected by an edge in the dual of the Johnson--Mehl graph. Note that since JMT is a planar graph, so is its dual, and thus $\mathcal D(I_0)$ has cardinality at most $3 \# I_0$ thanks to the Euler formula for planar graphs.

Now, let us define the distance of the closest point to the (space-time) origin in the Johnson--Mehl metric
\[ R'=\inf \{ r>0 \colon \exists i \in I \text{ with } d_{\rm J}((o,0),(X_i,T_i)) \leq r \}. \numberthis\label{R'def}\] 
Now, in the event $\{ R' \leq 1 \}$, we have $B_{1}^{\rm J}\cap \widetilde X \neq \emptyset$, and thus an application of Lemma~\ref{lem-JMindependence} for $a=b=1$ gives 
\[ \begin{aligned} W \leq  \# \{ \big((X_i,T_i);(X_j,T_j)\big) \in \Ecal \colon (X_i,T_i),(X_j,T_j) \in B_4^{\rm J} \} . \end{aligned} \] 
Thanks to \eqref{dualityargument}, the right-hand side is at most $ \# (\widetilde X \cap B_4^{\rm J})$. On the other hand, in the event $\{ R' > 1\}$, we can apply Lemma~\ref{lem-JMindependence} for $a=1$ and $b=R'$, which together with the convexity yields
\[ W \leq \# \{ \big((X_i,T_i);(X_j,T_j)\big) \in \Ecal \colon (X_i,T_i),(X_j,T_j) \in B_{R'+3}^{\rm J} \} . \numberthis\label{firstlargeR'} \]
Again, by stationarity of $X$ and absolute continuity of $\mu$, almost surely, we can further bound the right-hand side of \eqref{firstlargeR'} from above, which yields
\[ W \leq 3 \# (\widetilde X \cap B_{R'+3}) = 3+3\# (\widetilde X \cap (B_{R'+3} \setminus B_{R'})). \]
By assumption~\eqref{expmomentsvoidfiniteJM} we have $\E[R']<\infty$ and hence $\mathbb P(R'<\infty)=1$. 
We can thus estimate for all $\a>0$ using H\" older's inequality,
\begin{equation} \label{PoissonestimateJMT}
\begin{split}
\E[\exp(\a W)]\le\E\big[&\exp(3\a\# (\tilde X \cap B^{\rm J}_{4}))\big]\cr
&+\e^{3\a}\sum_{k\ge 1}\E\big[\exp\big(6\a\# (\tilde X \cap (B^{\rm J}_{k+4} \setminus B^{\rm J}_{k}))\big)\big]^{1/2}\P\big(\# (\tilde X \cap B^{\rm J}_{k})=0\big)^{1/2}.
\end{split}
\end{equation}
As above, the assumptions~\eqref{expmomentspointsfiniteJM} and~\eqref{expmomentsvoidfiniteJM} now guarantee summability. 
This proves Proposition~\ref{cor-number} for the number of edges and thus of Theorem~\ref{thm-expmoments} part (ii). 
\end{proof}

\subsubsection{Poisson--Delaunay tessellations: Proof of Theorem~\ref{thm-expmoments} part (iii) and Proposition~\ref{cor-number} part (i)}\label{sec-PDTproof}
The case of the DT is the most difficult one to handle, essentially since in this case, existence of points close to the origin does not automatically eliminate the influence of other distant points. To keep the argument simple, we thus only treat the case here where the underlying point process is a homogeneous PPP. 
Recall the definition of $W_a$ from~\eqref{Wadef}. Our first step towards the proof of Theorem~\ref{thm-expmoments} (iii) is to verify that there exists a fixed $\alpha>0$ such that $ \E[\exp(\alpha W_a)]<\infty $ holds for any $a>0$. Let us write $X^\l$ to indicate the intensity $\l$ in the underlying PPP and write $S^\l_{\rm D}=S_{\rm D}(X^\l)$ and $W^\lambda_a$ for the number of edges of $S^\l_{\rm D}$ intersecting with  $B_a$.
\begin{prop}\label{prop-PDTedges}
Let $a>0$ and $\lambda>0$. Then,
$\E[\exp(\alpha W^\lambda_a)]<\infty$
holds for all $\alpha<\frac{1}{6} \log \big( 1+ \frac{1}{72} \big)$.
\end{prop}
In particular, choosing $a=\lambda=1$, \eqref{expmomentsfinitenumberW} follows from this proposition for the Poisson--DT for small $\alpha>0$, which proves Proposition~\ref{cor-number} part (i) for the Poisson--DT. The proof rests on a comparison on the exponential scale. 
\begin{proof}
For $x \in \R^d$, let $Q_r(x)$ denote the box of side length $r$ centered at $x$. We define
\[ R=\min \{ r \in \N \colon r \geq 2a\text{ and } \forall z \in \Z^2 \text{ with } \Vert z \Vert_{\infty} = 2, Q_r(rz) \cap X^\lambda \neq \emptyset \}, \numberthis\label{NdefPDT}\]
the finest discretization of $\R^2$ into boxes such that every box in the $2$-annulus contains points. Note that $R$ is almost surely finite. For $k \in \N$ such that $k > \lceil 2a \rceil$,
\[ \begin{aligned}
\P(R \geq k) & \leq \P\big( \exists z \in \Z^2 \text{ with } \Vert z \Vert_{\infty} =2 \text{ such that } Q_{k-1}((k-1)z) \cap X^\lambda  = \emptyset \big) 
\\ & \leq \sum_{z \in \Z^2 \colon \Vert z \Vert=2} \P(  Q_{k-1}((k-1)z) \cap X^\lambda = \emptyset ) \leq 16 \P(  Q_{k-1}((k-1)\cdot(2,0)) \cap X^\lambda  = \emptyset ) 
\\ & \leq 16 \exp(-\lambda(k-1)^2).
\end{aligned} \numberthis\label{Restimate} \]
Note that once $k>\lceil 2a \rceil$, the right-hand side of \eqref{Restimate} does not depend on $a$. Since these terms are summable from $k=1$ to $\infty$, $\P(R=\lceil 2a \rceil)$ tends to one and thus $\E[R]$ tends to infinity as $a \to \infty$.

In the event $\{ R=k \}$ for some $k\ge 2a$, the points of $\partial Q_{3k}(o)$ are within a distance at most $\sqrt 2 k$ from the centroid of their Voronoi cell. Among these Voronoi cells, the neighboring ones are separated by a Voronoi edge and hence their cell centroids are Delaunay neighbors. The Delaunay edges connecting the centroids of the successive cells yield a closed path in the Delaunay graph surrounding $B_a$. This path defines a bounded region in which both endpoints of any Delaunay edge intersecting $B_a$ are located. Further, this region is fully contained in $Q_{3k}(o) \oplus B_{\sqrt 2k} \subset Q_{3k}(o) \oplus Q_{2\sqrt 2k}(o) \subset Q_{6k}(o)$. Hence, since the restriction of the Delaunay triangulation is a planar graph, using Euler's formula we arrive at
\[ W^\lambda _a \leq 3\# (X^\lambda  \cap Q_{6k}(o)). \]
Now we can use H\" older's inequality, the Laplace transform of a Poisson random variable and \eqref{Restimate} to estimate 
\[ 
\begin{aligned}
\E& [\exp(\alpha W^\lambda _a)]  \leq \sum_{k \geq 2a} \E \big[  \exp(3\alpha \# (X^\lambda  \cap Q_{6k}(o))) \mathds 1 \{ R = k \} \big] \\ & \leq  \sum_{k \geq 2a} \E \big[  \exp(6\alpha \# (X^\lambda  \cap Q_{6k}(o))) \big]^{1/2} \P( R = k )^{1/2} \leq
\sum_{k \geq \lceil 2a \rceil} \exp\big( 36 \lambda k^2 (\e^{6\alpha}-1) \big)\P( R = k )^{1/2}  \\ & \leq \exp\big( 144 \lambda (a+1)^2 (e^{6\alpha}-1) \big) + 4 \sum_{k=\lceil 2a \rceil+1}^\infty \exp\big( 36 \lambda k^2 (\e^{6\alpha}-1) \big) \exp\big( -\frac{1}{2}\lambda (k-1)^2 \big). 
%& \qquad +16\sum_{k= 4 \vee (\lceil 2a \rceil+1)}^{\infty} \exp\big( 36 \lambda k^2 (\e^{6\alpha}-1) \big) \exp\Big( -\frac{\lambda k^2}{2} \Big)
\end{aligned}
\]
But the right-hand side is finite for $\alpha<\frac{1}{6} \log \big( 1+\frac{1}{72} \big)$ for all $a>0$ and $\lambda>0$, as asserted. 
\end{proof}
We have the following corollary of Proposition~\ref{prop-PDTedges} for the total edge length.
\begin{cor}\label{cor-PDTlength}
Let $a>0$ and $\lambda>0$. Then for all $\alpha < \frac{1}{12a} \log \big( 1+\frac{1}{72} \big)$, $\E[\exp(\a|S^\lambda _{\rm D} \cap B_a|)]<\infty$.
\end{cor}
\begin{proof}
Since the edges of the Poisson--DT are straight line segments, any edge contributes to $|S^\lambda _{\rm D} \cap B_a|$ by at most $2a$. Hence,
\[ |S^\lambda _{\rm D} \cap B_a| \leq 2 a W^\lambda _a. \]
Now, $\E[\exp(\alpha (2a W^\lambda _a))]<\infty$ holds once $\E[\exp((2a\alpha)W^\lambda _a)]<\infty$. Thanks to Proposition~\ref{prop-PDTedges}, this holds as soon as $ \alpha < \frac{1}{12a} \log \big( 1+\frac{1}{72} \big)$, as wanted.
\end{proof}
Further, we have the following scaling relation for $\lambda,r>0$:
%\begin{lem}\label{lem-scaling}
%Let $\lambda,r>0$. Then we have the following identity in distribution
\[|S_{\rm D}^{\lambda} \cap B_1 |  =  \big| S_{\rm D}^{\l/r^2} \cap B_r \big|/r  \numberthis, \qquad\text{in distribution.} \label{scaleinvarianceDelaunay} \]
%\end{lem}
%\begin{proof}
Indeed, since $X^\lambda$, $X^{\lambda/r^2}$ are homogeneous PPPs with intensities $\lambda$, $\lambda/r^2$, respectively, we have that $X^{\lambda/r^2} \cap B_r$ equals $r(X^\lambda \cap B_1)$ in distribution. Thus, $S_{\rm D}^{\lambda/r^2} \cap B_r$ is equal to a rescaled version of $S_{\rm D}^\lambda \cap B_1$ in distribution where the length of each edge is multiplied by $r$. This implies the statement~\eqref{scaleinvarianceDelaunay}.
%\end{proof}
\begin{proof}[Proof of Theorem~\ref{thm-expmoments} part (iii)]
Let $\a, \lambda >0$. We need to show that $\mathbb E[\exp(\a|S_{\rm D}^\lambda\cap B_1|)]<\infty$. For this, let $\a_o<\frac{1}{12} \log \big( 1+\frac{1}{72} \big)$, then by Corollary~\ref{cor-PDTlength}, which works for arbitrary intensities, and the scale invariance~\eqref{scaleinvarianceDelaunay}, 
\begin{align*}
\infty&>\mathbb E[\exp(\a_o|S_{\rm D}^{n^2\lambda}\cap B_1|)]
=\mathbb E[\exp(\a_on^{-1}|S_{\rm D}^\lambda\cap B_n|)]=\mathbb E[\exp(\a_on^{1/2} n^{-3/2}|S_{\rm D}^\lambda\cap B_n|)],
\end{align*}
for all $n\in \N$. Now, for all $n\ge n_o$ where $n_o$ satisfies $\alpha_o n_o^{1/2}\ge \alpha$,
\begin{align*}
\infty&>\mathbb E[\exp(\a n^{-3/2}|S_{\rm D}^\lambda\cap B_n|)]\ge \mathbb E[\exp(\a n^{-3/2}|S_{\rm D}^\lambda\cap B_n|)\mathds 1\{n^{-3/2}|S_{\rm D}^\lambda\cap B_n|\ge |S^\lambda\cap B_1|\}],
\end{align*}
and hence 
\begin{align*}
\infty&>\mathbb E[\exp(\alpha |S_{\rm D}^\lambda\cap B_1|)\mathds 1\{n^{-3/2}|S_{\rm D}^\lambda\cap B_n|\ge|S_{\rm D}^\lambda\cap B_1|\}].
\end{align*}
Now, using Fatou's lemma, 
\begin{align*}
\infty&>\liminf_{n\uparrow\infty}\mathbb E[\exp(\alpha |S_{\rm D}^\lambda\cap B_1|)\mathds 1\{n^{-3/2}|S_{\rm D}^\lambda\cap B_n|\ge|S_{\rm D}^\lambda\cap B_1|\}]\\
&\ge \mathbb E[\exp(\alpha |S_{\rm D}^\lambda\cap B_1|)\liminf_{n\uparrow\infty}\mathds 1\{n^{-3/2}|S_{\rm D}^\lambda\cap B_n|\ge|S_{\rm D}^\lambda\cap B_1|\}]. 
\end{align*}
But, by ergodicity, $n^{-2}|S_{\rm D}^\lambda\cap B_n|\to \mathbb E[|S_{\rm D}^\lambda\cap B_1|]>0$ almost surely, and also almost surely, $n^{-1/2} |S_{\rm D}^\lambda\cap B_1|\to 0$, as $n$ tends to infinity. Hence, 
$$\liminf_{n\uparrow\infty}\mathds 1\{n^{-3/2}|S_{\rm D}^\lambda\cap B_n|\ge|S_{\rm D}^\lambda\cap B_1|\}=1,$$
also almost surely. But this implies the result. 
%Let us fix $\alpha$. Using \eqref{scaleinvarianceDelaunay}, it suffices to show that there exists $a>0$ such that
%\[ \mathbb E\Big[ \exp \Big( \frac{ \a}{a} \big| S_{\rm D}^{1/a^{2}} \cap B_{a} \big| \Big) \Big]<\infty \numberthis\label{rescaledwish} \]
%for some $a>0$.
%Thus, we only have to lift Corollary~\ref{cor-PDTlength} from sufficiently small $\a$ to all $\a$. 
%For $a>0$ let us define
%\[ \alpha_{\rm c}(a) = \frac{1}{12a} \log \big( 1 + \frac{1}{72} \big). \]
%Then, thanks to Corollary~\ref{cor-PDTlength}, $\E[\exp(\alpha|S_{\rm D} \cap B_1|)]<\infty$ for all $\alpha \in (0,\alpha_{\rm c}(1))$. 
%Let $r>0$ be sufficiently large such that $ \alpha/r<\alpha_{\mathrm c}(1)$. Note that $\alpha_{\rm c}(a) = \frac{1}{a} \alpha_{\rm c}(1)$. Further, observe that the value $\alpha_{\mathrm c}(1)$ is independent of the intensity parameter of the underlying PPP. These imply that for any $\lambda'>0$, we have
%\[ \mathbb E\Big[  \exp \Big( \frac{\a}{r} \big| S_{\rm D}^{\lambda} \cap B_r \big| \Big) \Big] < \infty. \]
%Choosing $\lambda=1/r^2$ implies \eqref{rescaledwish} with $a=r$ everywhere. This concludes the proof.
\end{proof}

\subsubsection{Line tessellations: Proof of Theorem~\ref{thm-expmoments} part (iv)}\label{sec-PLTproof}
We use the notation of Section~\ref{Exp_Setting}.
 Since for any line $l_i=\{ x \in \R^2 \colon x_1 \cos X_{i,2}+x_2 \sin X_{i,2}=X_{i,1} \}$ of $S_{\rm L}$ we have $|l_i \cap B_1| \leq 2$, it suffices to show that under the assumption \eqref{expmomentsfiniteLT} the number of lines of $S_{\rm L}$ intersecting with $B_1$ has exponential moments up to $\beta_\star$. Now, a line $l_i$ in $\R^2$ intersects with $B_1$ if and only if its distance parameter $X_{i,1}$ is at most one in absolute value, independently of its angle parameter $X_{i,2} \in [0,2\pi]$. By the assumption \eqref{expmomentsfiniteLT}, the number of such lines has exponential moments up to $\beta_\star$.
 \ProofEnde

\subsubsection{Manhattan grids: Proof of Theorem~\ref{thm-expmoments} for the MG}\label{sec-MGproof}
Since $B_1$ is a subset of $Q_1=Q_1(o)$, it suffices to verify the statement for $Q_1$ instead of $B_1$. Note that for any edge $E$ in $S_{\rm M}$, either $E \cap Q_1 = \emptyset$ or $| E \cap Q_1|=1$. Since $Y_{\rm v}$ and $Y_{\rm h}$ are independent, it follows that for all $\alpha>0$, we have
\begin{align*} \mathbb E[\exp(\alpha |S_{\rm M} \cap Q_1|)] & = \mathbb E[\exp\big(\alpha (\# (Y_{\mathrm v} \cap [-1/2,1/2]) +\#(Y_{\mathrm h} \cap [-1/2,1/2]))\big)] \\ & = \mathbb E[\exp\big(\alpha\# (Y_{\mathrm v} \cap [-1/2,1/2])\big)]\mathbb E[\exp\big(\alpha \big(\# Y_{\mathrm h} \cap [-1/2,1/2]\big)\big)].
\end{align*}
By assumption $\#(Y_{\mathrm h} \cap [-1/2,1/2])$ and $\#(Y_{\mathrm v} \cap [-1/2,1/2])$ have exponential moments, for $\a<\b_{\rm v}$, respectively $\a<\b_{\rm h}$, which implies exponential moments for $S_{\rm M}$ for $\a<\min\{\b_{\rm v},\b_{\rm h}\}$. 
 \ProofEnde

\subsubsection{Number of cells: Proof of Proposition~\ref{cor-number} part (ii)}\label{sec-cells}
\begin{proof}[Proof of Proposition~\ref{cor-number} part (ii)]
Note that any edge of the VT, DT or JMT that intersects with $B_1$ is adjacent to precisely two cells intersecting with $B_1$, whereas if $W=0$, then $V=1$, and thus we have the trivial bound $ V \leq 2 W+1$. Thus, the assertion \eqref{expmomentsfinitenumberV} for any given $\alpha/2>0$ follows from the assertion \eqref{expmomentsfinitenumberW} for the same $\alpha$.
\end{proof}
\subsection{Nested tessellations: Proof of Corollary~\ref{cor-iterated} and Proposition~\ref{prop-iterated-MG}}
\begin{proof}[Proof of Corollary~\ref{cor-iterated}.]
We write $S'$ for a fixed tessellation process that equals $S_i$, $i \in J$, in distribution, and we define $V$ according to \eqref{Vstardef} for the first-layer tessellation $S_o$, so with the index set $J$ being such that $S_o$ has cells $(C_i)_{i \in J}$. For $\alpha,\beta>0$, let us write
\[ M_\alpha=\E[\exp(\alpha |S' \cap B_1|)] \qquad \mbox{\text{and}} \qquad N_\beta = \E[\exp(\beta V)], \]
where $M_\alpha, N_\beta$ are defined as elements of $[0,\infty]$. Then, we need to show \eqref{first-cor-iterated} that if $M_\alpha<\infty$ and $N_\beta<\infty$ for all $\alpha,\beta>0$, then  $\E[\exp(\gamma |S_\mathrm N \cap B_1|)]<\infty$ holds for all $\gamma>0$, and \eqref{second-cor-iterated} if there exists $\alpha,\beta>0$ such that $M_\alpha<\infty$ and $N_\beta<\infty$, then there exists $\gamma>0$ such that $\E[\exp(\gamma |S_\mathrm N \cap B_1|)]<\infty$. First, using H\" older's inequality, we can separate the first from the second layer process, 
\begin{align*}
 \E & [\exp(\alpha|S_{\rm N} \cap B_1|)]  \leq \E \Big[ \exp \Big( 2\alpha \sum_{\begin{smallmatrix} i \in J \colon C_i \cap B_1 \neq \emptyset \end{smallmatrix}} |S_i \cap C_i \cap B_1| \Big) \Big]^{\frac{1}{2}}  \E [\exp(2\alpha|S_o \cap B_1|)]^{\frac{1}{2}}. 
\end{align*}
%where by assumption $\E [\exp(2\alpha|S_o \cap B_1|)]<\infty$. \color{red} It is not clear which $\alpha$ is meant here because there are two different cases. Maybe we can delete the statement that the expectation is finite, it will soon be said once again. \color{black} 
For the first factor on the right-hand side, note that we can bound
\begin{align*}
&\E \Big[ \exp \Big( 2\alpha \sum_{\begin{smallmatrix} i \in J \colon C_i \cap B_1 \neq \emptyset \end{smallmatrix}} |S_i \cap C_i \cap B_1| \Big) \Big] 
 \\ & = \E \Big[ \E \Big[ \exp \Big( 2\alpha\sum_{\begin{smallmatrix} i \in J \colon  C_i \cap B_1 \neq \emptyset \end{smallmatrix}} |S_i \cap C_i \cap B_1| \Big) \Big| S_o \Big]\Big]
 \leq \E \Big[ \E \Big[ \exp \Big( 2\alpha \sum_{i \in J \colon C_i \cap B_1 \neq \emptyset} |S_i \cap B_1| \Big) \Big| S_o \Big]\Big] \\ & = \E \Big[ \prod_{i \in J \colon C_i \cap B_1 \neq \emptyset} \E \Big[ \exp \Big( 2\alpha  |S_i \cap B_1| \Big) \Big| S_o \Big]\Big]  = \E \big[ M_{2\alpha}^{V} \big]=\E \big[ \exp(V\log M_{2\alpha}) \big] = N_{\log M_{2\alpha}},
\end{align*}
as an inequality in $[0,\infty]$. From this, \eqref{first-cor-iterated} follows immediately. As for \eqref{second-cor-iterated}, let us assume that $M_\alpha<\infty$ holds for some $\alpha>0$ and $N_\beta <\infty$ holds for some $\beta>0$. Then, the moment generating function $\R \to [0,\infty]$, $\beta \mapsto N_\beta$ is continuous (in fact, infinitely many times differentiable) in an open neighborhood of $0$, which implies that $\lim_{\beta \to 0} N_{\beta}=N_0=1$. Analogous arguments imply that $\lim_{\alpha \to 0} \log M_{\alpha}=0$. Hence, there exists $\alpha>0$ such that $N_{\log M_\alpha}<\infty$, which implies \eqref{second-cor-iterated}.
\end{proof}
\begin{proof}[Proof of Proposition~\ref{prop-iterated-MG}.]
We verify the statement with $B_1$ replaced by $Q_1$ in \eqref{expmomentsfinite}, which suffices thanks to the fact that $B_1 \subset Q_1$. According to the assumptions of the proposition, let the first-layer tessellation $S_o$ be a MG satisfying \eqref{expmomentsfinite} for all $\alpha>0$,  and let us write $Y^o=(Y_\mathrm v^o, Y_\mathrm h^o)$ for the corresponding pair of point processes on $\R$. We can enumerate the points of $Y_{\rm v}^o \cap [-1/2,1/2]$ in increasing order as $Y_{\rm v}^o \cap [-1/2,1/2]=(P^i)_{i=1}^{N_{\rm v}}$.
% where $N_{\rm v}$ is a Poisson random variable with parameter $\lambda_{\rm v}$, and given $N_{\rm v}$, $(P^i)_{i=1}^{N_{\rm v}}$ are the increasing order statistics of $N_{\rm v}$ i.i.d.~uniformly distributed random variables in $[-1/2,1/2]$. 
Similarly, we can enumerate  the points of $Y_{\rm h}^o \cap [-1/2,1/2]$ in increasing order as $Y_{\rm h}^o \cap [-1/2,1/2]=(Q^j)_{j=1}^{N_{\rm h}}$.
%, where $N_{\rm h}$ is a Poisson random variable with mean $\lambda_{\rm h}$. 
We further write $P^0=Q^0=-1/2$ and $P^{N_{\rm v}+1}=Q^{N_{\rm h}+1}=1/2$. Note that $\sum_{i=1}^{N_{\mathrm v}+1} (P^i-P^{i-1}) = \sum_{j=1}^{N_{\mathrm h}+1} (Q^j-Q^{j-1}) = 1$.

Now, the collection of cells of $S_o$ intersecting $Q_1$ is given as $(C_{i,j})_{i=1,\ldots,N_{\rm v}+1, j=1,\ldots,N_{\rm h}+1}$, where $C_{i,j}$ is the open rectangle $(P^{i-1},P^{i}) \times (Q^{j-1},Q^{j})$. We write $S_{i,j}$ for the second-layer tessellation corresponding to $S_{\rm N}$ in the cell $C_{i,j}$ and $Y^{i,j}=(Y^{i,j}_{\rm v},Y^{i,j}_{\rm h})$ for the associated pair of Poisson processes on $\R$. Here, there exist $\lambda_{\rm v},\lambda_{\rm h}>0$ such that for all $i \in \{ 1,\ldots,N_{\rm v}+1 \}$ and for all $j\in \{1,\ldots,N_{\rm h}+1 \}$, $Y^{i,j}_{\rm v}$ has intensity $\lambda_{\rm v}$ and $Y^{i,j}_{\rm h}$ has intensity $\lambda_{\rm h}$. Now note that for all $i \in \{ 1,\ldots,N_{\rm v}+1 \}$ and for all $j\in \{1,\ldots,N_{\rm h}+1 \}$, all vertical edges of $S_{i,j}$ intersect $C_{i,j}$ in a segment of length $P^i-P^{i-1}$ and all horizontal edges of $S_{i,j}$ intersect $C_{i,j}$ in a segment of length $Q^j-Q^{j-1}$. Thus, we obtain that
\[
\begin{aligned}   |S_{\rm N} \cap Q_1 | =  |S_o \cap Q_1 | +& \sum_{i=1}^{N_{\rm h}+1} (P^i-P^{i-1})  \sum_{j=1}^{N_{\rm v} + 1} \#\big(Y^{i,j}_{\rm v} \cap (Q^{j-1},Q^j)\big) \\  &+ \sum_{j=1}^{N_{\rm v} + 1} (Q^j-Q^{j-1}) \sum_{i=1}^{N_{\rm h}+1} \# \big(Y^{i,j}_{\rm h} \cap (P^{i-1},P^i)\big).
\end{aligned} \]
By Hölder's inequality, it suffices to verify the existence of all exponential moments for each of the three terms on the right-hand side separately. The first term has all exponential moments thanks to the assumption of Proposition~\ref{prop-iterated-MG}. Further, by symmetry between the second and the third term, it suffices to show existence of all exponential moments for one of them; we will consider the second term.

Since for fixed $i \in \{1,\ldots,N_{\mathrm h} + 1 \}$, $\# (Y_v^{i,j} \times (Q^{j-1},Q^j))_{j=1,\ldots,N_{\mathrm v}+1}$ are independent Poisson random variables with parameters summing up to $\lambda_{\rm v}$, it follows that their superposition $N_i = \sum_{j=1}^{N_{\mathrm v}+1} \# (Y_{\mathrm v}^{i,j} \cap (Q^{j-1},Q^j)$ is a Poisson random variable with parameter $\lambda_{\rm v}$. Further, conditional on $(P^i)_{i=1}^{N_{\mathrm h}}$, $(N_i)_{i=1}^{N_{\mathrm h}+1}$ are independent.

Now, fix $\alpha>0$, and let $K_\alpha>0$ be such that for all $x \in (-\infty,\alpha]$ we have $\exp(x)-1 \leq K_\alpha x$. Using that $P^i-P^{i-1} \leq 1$ for all $i$ and $\sum_{i=1}^{N_{\rm h}+1} (P^i-P^{i-1})=1$, we estimate
\begin{align*}
\E & \Big[ \exp \Big( \alpha \sum_{i=1}^{N_{\rm h}+1} (P^i-P^{i-1})  \sum_{j=1}^{N_{\rm v} + 1} \# (Y^{i,j}_{\rm v} \cap (Q^{j-1},Q^j)) \Big)\Big]  =  \E\Big[ \exp \Big( \alpha \sum_{i=1}^{N_{\rm h}+1} (P^i-P^{i-1}) N_i  \Big) \Big] 
 \\ & =  \E \Big[ \E\Big[ \exp \Big( \alpha \sum_{i=1}^{N_{\rm h}+1} (P^i-P^{i-1}) N_i  \Big) \Big| (P^i)_{i=1}^{N_{\rm h}} \Big] \Big] = \mathbb E \Big[ \prod_{i=1}^{N_{\mathrm h} + 1 } \E\Big[ \exp \Big( \alpha (P^i-P^{i-1}) N_i  \Big) \Big| (P^i)_{i=1}^{N_{\rm h}} \Big] \Big] 
\\  &=  \mathbb E \Big[ \prod_{i=1}^{N_{\mathrm h} + 1 } \exp \big( \lambda_{\rm v} ( \exp( \alpha (P^i-P^{i-1}))-1)   \big)  \Big]  
 \leq  \mathbb E \Big[ \prod_{i=1}^{N_{\mathrm h} + 1 } \exp \big( K_\alpha  \lambda_{\rm v} \alpha (P^i-P^{i-1})  \big)  \Big] \\ &= \mathbb E \Big[ \exp \Big(  \sum_{i=1}^{N_{\mathrm h} + 1 } K_\alpha  \lambda_{\rm v} \alpha (P^i-P^{i-1})  \Big)  \Big] =  \exp \big( K_\alpha  \lambda_{\rm v} \alpha  \big). 
\end{align*}
Since the right-hand side is finite, we conclude the proof of the proposition.
\end{proof}

\subsection{Palm versions of tessellations: Proof of Corollary~\ref{cor-Palm}}\label{sec-Palm}
We handle each case separately.
\begin{proof}[Proof of Corollary~\ref{cor-Palm} for the Poisson--VT]
Corollary~\ref{cor-Palm} follows directly from Lemma~\ref{lem-conditioningonapoint} and the Slivnyak--Mecke theorem (see e.g.~\cite[Section 9.2]{LP17}). Indeed, since Lemmas~\ref{lem-conditioningonapoint} uses no information about the distribution of $X$ but only the definition of a Voronoi tessellation, these lemmas remain true after replacing $S^*$ by $S$. Next, the Palm version $X^*$ of the underlying PPP equals $X \cup \{ o \}$ in distribution by the Slivnyak--Mecke theorem, in particular, it contains $o$ almost surely. Thus, using the aforementioned versions of Lemma~\ref{lem-conditioningonapoint} (for $a=b=1$), we deduce that $|S_{\rm V} \cap B_1|$ is stochastically dominated by $ 2\pi (\# (X \cap B_4)+1)$. This random variable has all exponential moments, hence the corollary.
\end{proof}
\begin{proof}[Proof of Corollary~\ref{cor-Palm} for the Poisson--JMT]
This is analogous to the proof for the Poisson--VT where instead of Lemma~\ref{lem-conditioningonapoint} we use the Lemma~\ref{lem-JMindependence}.
\end{proof}
\begin{proof}[Proof of Corollary~\ref{cor-Palm} for the Poisson--DT]
Note that the random radius $R$ defined in \eqref{NdefPDT} is invariant under changing $X$ to $X \cap \{ o \}$ in its definition. Hence, using the Slivnyak--Mecke theorem, one can first verify Proposition~\ref{prop-PDTedges} with $W_a$ replaced by the number of edges of $S_{\rm D}^*$ intersecting with $B_a$, then one can prove that Corollary~\ref{cor-PDTlength} holds with $S_{\rm D}$ replaced by $S_{\rm D}^*$ and Lemma~\ref{lem-scaling} holds with $S_{\rm D}^\lambda$ replaced by its Palm version $(S_{\rm D}^\lambda)^*$ for all $\lambda>0$, and then one can complete the proof of Corollary~\ref{cor-Palm} for the Poisson--DT analogously to the final part of the proof of Theorem~\ref{thm-expmoments} (iii).
\end{proof}
\begin{proof}[Proof of Corollary~\ref{cor-Palm} for the Poisson--LT] Recall from Section~\ref{Exp_Setting} that $S^*$ equals $S_L(X^{*})$ where $X^{*} = X \cup \{ (0,\Phi) \}$, with $\Phi$ being a uniform random angle in $[0,\pi)$ that is independent of $X$.Thus, $S^* = S \cup \{ l \}$, where
$ l=\{ x \in \R^2 \colon x_1 \cos \Phi + x_2 \sin \Phi = 0 \}. $
Since the intersection of $l$ with $B_1$ has length 2, the corollary in the case of a Poisson--LT follows directly from Theorem~\ref{thm-expmoments} part (iv). 
\end{proof}
\begin{proof}[Proof of Corollary~\ref{cor-Palm} for the MG]
We verify the statement with $B_1$ replaced by its superset $Q_2$. First, let us write $Y_\mathrm v^*$ and $Y_\mathrm h^*$ for the Palm versions of $Y_\mathrm v$ and $Y_\mathrm h$.  Here, $Y_{\mathrm v}^*$ is defined via the property \cite[Section 2.2]{HJC18} that
\[ \E \big[ f(Y_{\rm v}^*)  \big] =\frac{1}{\lambda_{\mathrm v}}  \E \Big[ \sum_{X_i \in Y_{\rm v} \cap [0,1]} f(Y_{\rm v}-X_i) \Big] \]
for any measurable function $f$ on the space of $\sigma$-finite counting measures on $\R$ to $[0,\infty)$.
Then the Palm version $S_{\rm M}^*$ is given according to \eqref{PalmMG}. It suffices to verify that $Y_{\mathrm v}^* \times [-1,1]$ and $Y_{\mathrm h}^* \times [-1,1]$ have all exponential moments. Indeed, using this and the mutual independence of $Y_{\mathrm v}$, $Y_{\mathrm h}$, and $U$, the proof of Corollary~\ref{cor-Palm} for the MG can be completed analogously to the proof of Theorem~\ref{thm-expmoments} part (v) in Section~\ref{sec-MGproof}. We only consider $Y_{\mathrm v}^*$, the proof for $Y_{\mathrm h}^*$ is analogous.
For $\alpha>0$ we have
\begin{align*}
& \E[\exp\big(\alpha \# (Y_{\mathrm v}^* \cap [-1,1])\big)] =  \frac{1}{\lambda_{\rm v}}\E \Big[ \sum_{X_i \in Y_\mathrm v \cap [0,1]} \exp\big(\alpha \# ((Y_\mathrm v-X_i) \cap [-1,1])\big) \Big] \\
& =   \frac{1}{\lambda_{\rm v}}\E \Big[ \sum_{X_i \in Y_\mathrm v \cap [0,1]} \exp\big(\alpha \# (Y_\mathrm v \cap [X_i-1,X_i+1])\big) \Big] \\\ & \leq  \frac{1}{\lambda_{\rm v}} \E \Big[ \sum_{X_i \in Y_\mathrm v \cap [0,1]} \exp\big(\alpha \# (Y_\mathrm v \cap [-2,2])\big) \Big] 
 = \frac{1}{\lambda_{\rm v}} \E\big[\# (Y_{\mathrm v} \cap [0,1])  \exp\big(\alpha \# (Y_\mathrm v \cap [-2,2])\big)\big] \\ & \leq \frac{1}{\lambda_{\rm v}} \E[\# (Y_{\mathrm v} \cap [0,1])^2]^{1/2} \E[ \exp\big(2\alpha \# (Y_\mathrm v \cap [-2,2])\big)]^{1/2} <\infty,
\end{align*}
where in the first inequality of the last line we used H\" older's inequality.
\end{proof}

\medskip
As above, note that weaker assumptions on the exponential moments of $Y_{\mathrm v},Y_{\mathrm h}$ imply lower exponential moments for $S^*_{\rm MG}$. 

%\section{Discussion}\label{Exp_Discussion}
%\input{Exp_Discussion}

\subsection*{Acknowledgements}
The authors thank A.~Hinsen, C.~Hirsch and W.~K\" onig for interesting discussions and comments. Moreover, we thank an anonymous reviewer for suggestions regarding simplifications of the proof of Theorem~\ref{thm-expmoments} (i)-(iii) as well as possible generalizations, which helped to substantially improve the manuscript. Further, we thank two more anonymous reviewers for insightful suggestions that lead to further improvements. The first author was supported by the Deutsche Forschungsgemeinschaft (DFG, German Research Foundation) under Germany's Excellence Strategy – MATH+ : The Berlin Mathematics Research Center, EXC-2046/1 – project ID: 390685689. The second author was supported by a Phase II scholarship from the Berlin Mathematical School.

\end{document}